\documentclass[12pt]{amsart}

\usepackage[dvips]{graphicx}
	\usepackage{amsmath,amssymb,amsthm,wrapfig,amsfonts,enumerate,latexsym}

	\newtheorem{dfn}{Definition}[section]
	\newtheorem{thm}[dfn]{Theorem}
	\newtheorem{prop}[dfn]{Proposition}
	\newtheorem{cor}[dfn]{Corollary}
	\newtheorem{lem}[dfn]{Lemma}
	\newtheorem{rem}[dfn]{Remark}
	\newtheorem{ex}[dfn]{Example}
 	\newtheorem{claim}[dfn]{Claim}

	\newtheorem{ack}{Acknowledgements\!\!}

	\newcounter{yon}
	\numberwithin{equation}{section}

	\def\notin{\not\in}

	\newcommand{\dist}{\mathop{\mathit{d}} \nolimits}
	
	\newcommand{\diam}{\mathop{\mathrm{diam}} \nolimits}

	\newcommand{\sep}{\mathop{\mathrm{Sep}} \nolimits}
	\newcommand{\ler}{\mathop{\mathrm{LeRad}} \nolimits}

	\newcommand{\me}{\mathop{\mathrm{me}}      \nolimits}
	\newcommand{\supp}{\mathop{\mathrm{Supp}}    \nolimits}
	\newcommand{\rlip}{\mathop{\mathrm{Lip}}     \nolimits}
	\newcommand{\oblip}{\mathop{\stackrel{\rlip_1}{\longrightarrow}}    \nolimits}

	\newcommand{\obsc}{\mathop{\mathrm{ObsCRad}}               \nolimits}
	\newcommand{\crad}{\mathop{\mathrm{CRad}}           \nolimits}

	\newcommand{\obin}{\mathop{\mathrm{Obs}L^p\mathrm{\text{-}Var}}
	\nolimits}
	\newcommand{\obinin}{\mathop{\mathrm{Obs}L^2\mathrm{\text{-}Var}}
	\nolimits}

	\newcommand{\ric}{\mathop{\mathit{Ric}}        \nolimits}

    \newcommand{\gcon}[1]{\mathop{C_{#1, \mathrm{Ga}}}
    \nolimits}
    \newcommand{\econ}[1]{\mathop{C_{#1, \mathrm{exp}}}
    \nolimits}
     \newcommand{\gconn}[1]{\mathop{c_{#1, \mathrm{Ga}}}
    \nolimits}
    \newcommand{\econn}[1]{\mathop{c_{#1, \mathrm{exp}}}
    \nolimits}

\begin{document}

	\title[]
	{Observable concentration of mm-spaces into nonpositively
	curved manifolds}
	\author[Kei Funano]{Kei Funano}
	\address{Mathematical Institute, Tohoku University, Sendai 980-8578, JAPAN}
	\email{sa4m23@math.tohoku.ac.jp}
\thanks{This work was partially supported by
Research Fellowships of the Japan Society for the Promotion of Science for Young Scientists.}
	\dedicatory{}
	\date{\today}

	\maketitle


	\setlength{\baselineskip}{5mm}

	\begin{abstract}The measure concentration property of an mm-space $X$ is roughly described as that any $1$-Lipschitz map
	 on $X$ to a metric space $Y$ is almost close to a constant map. The target space $Y$ is called the screen. The case of
	 $Y=\mathbb{R}$ is widely studied in many literature (see
	 \cite{gromov}, \cite{ledoux}, \cite{mil2}, \cite{milsch},
	 \cite{sch}, \cite{tal}, \cite{tal2} and their
	 references). 
	 M. Gromov developed the theory of measure concentration in
	 the case where the screen $Y$ is not necessarily $\mathbb{R}$
	 (cf.~\cite{gromovcat}, \cite{gromov2}, \cite{gromov}). In this
	 paper, we consider the case where the screen $Y$ is a nonpositively curved manifold. We also show that if the screen $Y$ is so
	 big, then the mm-space $X$ does not concentrate.
	\end{abstract}

\section{Introduction}
	Let $\mu_{n}$ be the volume measure on the $n$-dimensional unit sphere $\mathbb{S}^n$ in $\mathbb{R}^{n+1}$ normalized as $\mu_n (\mathbb{S}^n)=1$. In $1919$, P. L\'{e}vy proved that
	for any $1$-Lipschitz function $f:\mathbb{S}^n \to \mathbb{R}$ and any $\varepsilon>0$, the inequality
	\begin{align}\label{11922960}
\mu_n\big(\{ x\in \mathbb{S}^n \mid |f(x)-m_f|\geq \varepsilon    \}\big)\leq 2\ e^{-(n-1){\varepsilon}^2/2} 
	 \end{align}holds, where $m_f$ is some constant determined by $f$. For any fixed $\varepsilon>0$ the right-hand side of the above
	inequality converges to zero as $n\to \infty$. This
	means that any $1$-Lipschitz function on $\mathbb{S}^n$ is
	almost close to a constant function for sufficiently large $n\in
	\mathbb{N}$. This high dimensional concentration phenomenon of
	functions was
	first extensively used and emphasized by V. D. Milman in his
	investigation of asymptotic geometric analysis. He used
	L\'{e}vy's inequality (\ref{11922960}) for a short proof of
	Dvoretzky's theorem on Euclidean section of convex bodies (cf.~\cite{mil}). The same year later he extended L\'{e}vy's
	result to some nonspherical spaces and then pushed forward the
	idea of concentration of functions as a general unifying principle
	(cf.~\cite{mill}, \cite{mil1}). Nowadays, the concentration theory of
	functions is widely studied in many literature and blend with
	various areas of mathematics (see \cite{gromov}, \cite{ledoux}, \cite{mil2},
	\cite{milsch}, \cite{sch}, \cite{tal}, \cite{tal2} and references therein for
	further information).

In
	1999, M. Gromov 
	established a theory of concentration of maps 
	into general metric spaces by introducing the notion of the observable diameter in \cite{gromov}. He settled the following definition. 

	\begin{dfn}\upshape
	 Let $Y$ be a metric space and $\nu_Y$ a Borel measure on $Y$ such that $m:=\nu_Y(Y)<+\infty$. 
	 We define for any $\kappa >0$
	 \begin{align*}
	  \diam (\nu_Y , m-\kappa):= \inf \{ \diam Y_0 \mid Y_0 \subseteq Y \text{ is a Borel subset such that }\nu_Y(Y_0)\geq m-\kappa\}
	  \end{align*}and call it the \emph{partial diameter} of $\nu_Y$. 
	 \end{dfn}
	An \emph{mm-space} is a triple $(X,\dist,\mu)$, where $\dist$ is a complete separable metric on
	a set $X$ and $\mu$ a finite Borel measure on $(X,\dist)$.

	\begin{dfn}[Observable diameter]\upshape Let $(X,\dist,\mu)$ be an mm-space and $Y$ a metric space. For any $\kappa >0$ we
	 define the \emph{observable diameter} of $X$ by 
	 \begin{align*}
	  \diam (X\oblip Y, m-\kappa):=
	   \sup \{ \diam (f_{\ast}(\mu),m-\kappa) \mid f:X\to Y \text{ is a }1 \text{{\rm -Lipschitz map}}  \},
	  \end{align*}where $f_{\ast}(\mu)$ stands for the push-forward measure of $\mu$ by $f$. The target metric space $Y$ is called the \emph{screen}. 
	 \end{dfn}
	The idea of the observable diameter comes from the quantum and statistical
	mechanics, that is, we think of $\mu$ as a state on a configuration space $X$ and $f$ is interpreted as an observable. We define a sequence $\{  X_n   \}_{n=1}^{\infty}$ of
	mm-spaces is a \emph{L\'{e}vy family} if $\diam (X_n\oblip \mathbb{R},m_n-\kappa)\to 0$ as $n\to \infty$ for any $\kappa >0$, where
	$m_n$ is the total measure of the mm-space $X_n$. This is equivalent to that for any
	$\varepsilon >0$ and any sequence $\{ f_n :X_n \to \mathbb{R}\}_{n=1}^{\infty}$ of $1$-Lipschitz functions, we have 
	\begin{align*}
	 \mu_n(\{ x\in X_n \mid |f_n(x) - m_{f_n}|\geq \varepsilon
	 \})\to 0 \text{ as } n\to \infty,
	\end{align*}where $m_{f_n}$ is a some constant determined by
	 $f_n$. The L\'{e}vy families are first introduced and analyzed
	from a geometric point of view by Gromov and Milman in
	\cite{milgro}. The inequality (\ref{11922960}) shows that the sequence $\{ \mathbb{S}^n\}_{n=1}^{\infty}$
	is a L\'{e}vy family. Gromov proved in \cite{gromov} that $\diam (X_n \oblip Y, m_n - \kappa) \to 0$ as $n\to \infty$ for any
	$\kappa >0$ and for a L\'{e}vy family $\{ X_n\}_{n=1}^{\infty}$
	if the screen $Y$ is a compact metric space or a Euclidean
	space.  He also discussed the case where the screens are
	Euclidean spaces whose dimensions go to infinity by considering the
	barycenters of the push-forward measures (see Theorem \ref{g}). His barycenter method also
	goes well if the screens are nonpositively curved manifolds (see Section $4$). In \cite{gromov2}, he considers and analyzes the questions of isoperimetry of waists and concentration
	for maps from a unit sphere to a Euclidean space (see Theorem \ref{isowa}). In the recent work \cite{ledole}, M. Ledoux and K. Oleszkiewicz estimated
	the observable diameter $\diam (X \oblip \mathbb{R}^k,m-\kappa)$
	provided that the mm-space $X$ has a Gaussian concentration (see
	Theorem \ref{thledole}). In our
	previous paper \cite{funano}, the author studied the case where the screen $Y$ is a metric space
	with a doubling measure. 

	In this paper, inspired by Gromov's study, we
	study concentration phenomenon of maps into nonpositively curved
	manifolds. In
	particular, we consider the case where dimensions of screens
	go to infinity. We denote by $\mathcal{NM}^n$ the set of all
	$n$-dimensional complete simply connected Riemannian manifolds
	with nonpositive sectional curvature. One of our main theorems is the following:
	\begin{thm}\label{kuzuninaritaiwan}Let $\{ X_n \}_{n=1}^{\infty}$ be a sequence of
	 mm-spaces. We assume that a sequence $\{
	 a(n)\}_{n=1}^{\infty}$ of natural numbers satisfies that for
	 any $\kappa>0$
	 \begin{align}\label{denwanidenwa}
	  \diam (X_n \oblip \mathbb{R}^{a(n)} ,m_n-\kappa) \to 0 \text{
	  as }n\to \infty.
	  \end{align}Then for any $\kappa>0$, we have
	 \begin{align*}
	  \sup \{    \diam (X_n \oblip N,m_n-\kappa) \mid N\in
	  \mathcal{NM}^{a(n)}       \}\to 0 \text{ as }n\to \infty.
	  \end{align*}
	 \end{thm}
	In the proof of Theorem \ref{kuzuninaritaiwan}, we find a point in a screen
	$N$ which is a kind
	of barycenter of the push-forward measure on $N$,
	and prove that the measure concentrates to the point by the
	delicate discussions comparing $N$ with both an Euclidean space and a real hyperbolic space.

In \cite{gromovcat}, Gromov introduced
	the notion of $L^p$-concentration of maps from mm-spaces into
	metric spaces. We recall its definition. Let $(X,\dist_X,\mu_X)$ be an mm-space and $(Y,\dist_Y)$ a metric
space. Given a Borel measurable map $f:X\to Y$ and $p\in (0,+\infty)$, we put 
\begin{align*}
&V_p(f):=  \Big(\int \int_{X\times X}
 \dist_Y\big(f(x),f(x')\big)^p   \ d\mu_X(x)\ d \mu_X(x')\Big)^{1/p},\\
&V_{\infty}(f):=   \inf \big\{  r>0 \mid
 (\mu_X \times \mu_X)( \{  (x,x')\in X\times X   \mid
 \dist_Y(f(x),f(x'))\geq r  \}) =0 \big\}.
\end{align*}

Let $ \{  X_n        \}_{n=1}^{\infty}$ be a
 sequence of mm-spaces and $\{  Y_n   \}_{n=1}^{\infty}$ a sequence of
 metric spaces. For any $p\in (0,+\infty]$, we say that a sequence $\{ f_n :X_n \to Y_n
 \}_{n=1}^{\infty}$ of Borel measurable maps \emph{$L^p$-concentrates}
 if $V_{p}(f_n)\to 0$ as $n\to \infty$. 

 We easily see that $L^p$-concentration of maps implies the
	concentration of maps (see Lemma \ref{var1} and Corollary \ref{hayakuowaritai1}). In
	\cite{gromovcat}, Gromov got several $L^2$-concentration
	inequalities of maps (see for example, Theorem \ref{sugeyo}).

    Given an mm-space $X$ and a metric space $Y$ we define
\begin{align*}
\obin_Y(X):= \sup \{   V_p(f) \mid f:X\to Y \text{ is a }1
 \text{-Lipschitz map}                  \},
\end{align*}and call it the \emph{observable} $L^p$\emph{-variation} of
 $X$. We get the following $L^p$-concentration result:

 \begin{thm}\label{koukaithm}Let $X$ be an mm-space and $p\geq 1$. Then,
  we have
  \begin{align}\label{katouta}
   \obin_N(X)\leq 2\obin_{\mathbb{R}^n}(X)
   \end{align}for any $n\in \mathbb{N}$ and $N\in \mathcal{NM}^n$. In
  particular, if a sequence $\{X_n\}_{n=1}^{\infty}$ of mm-spaces and a
  sequence $\{a(n)\}_{n=1}^{\infty}$ of natural numbers satisfy that
  \begin{align}\label{warauyokogao}
   \obin_{\mathbb{R}^{a(n)}}(X_n)\to 0 \text{ as }n\to \infty
   \end{align}for some $p\geq 1$, then we have
  \begin{align}\label{warauyokogao1}
   \sup \{  \obin_{N}(X_n) \mid N\in \mathcal{NM}^{a(n)}\}\to 0 \text{
   as }n\to \infty.
   \end{align}
  \end{thm}

In \cite[Section 13]{gromovcat}, Gromov got $\obinin_N(X)\leq \sqrt{
2}\obinin_{\mathbb{R}^n}(X)$. Note that this Gromov's inequality implies
better estimate than our inequality (\ref{katouta}) (see also Remark
\ref{katoukatou}). Our proof is an analogue to his proof.

The following theorem includes the case of $0<p<1$:
\begin{thm}\label{kameda1}Let $\{ X_n\}_{n=1}^{\infty}$ be a sequence of
 mm-spaces with finite diameter. Assume that a sequence $\{  a(n)\}_{n=1}^{\infty}$ of
 natural numbers satisfies that $\sup\limits_{n\in \mathbb{N}}m_n <+\infty$ and \begin{align}\label{hayakuowaritai2}
\diam \Big( X_n \oblip \mathbb{R}^{a(n)},m_n -\frac{\kappa}{(\diam
 X_n)^{p}}\Big) \to 0 \text{ as }n\to \infty
\end{align}for any $\kappa >0$ and some $p>0$. Then we have
\begin{align}\label{hayakukaeritai4}
\sup \{ \obin_{N}(X_n) \mid N\in \mathcal{NM}^{a(n)}  \} \to 0 \text{ as
 }n\to \infty.
\end{align} 
\end{thm}
	Note that the assumption (\ref{hayakuowaritai2}) directly
	implies (\ref{denwanidenwa}). Actually, (\ref{hayakuowaritai2})
	leads to (\ref{warauyokogao}) (see Corollary \ref{waruiyokan}). We do not know
	whether the assumption (\ref{warauyokogao}) for $0<p<1$ implies (\ref{hayakukaeritai4})
	or not. By using Gromov's observation in \cite{gromov}, we prove a
	similar result of Theorems \ref{koukaithm} and \ref{kameda1} (see
	Proposition \ref{shishi} and Remark \ref{redo1}).

The next proposition says that if the screen $Y$ is so big that the
	mm-space $X$ which satisfies some homogeneity condition, can isometrically be embedded into $Y$, then its observable diameter $\diam(X\oblip Y,m-\kappa)$ is not close to zero.
	\begin{prop}\label{kantan}
	 Let $\{  (X_n,\dist_n ,\mu_n)\}_{n=1}^{\infty}$ be a L\'{e}vy
	 family with $\inf\limits_{n\in \mathbb{N}}m_n>0$ and $\{ Y_n \}_{n=1}^{\infty}$ a sequence of metric spaces satisfying the
	 following conditions $(1)-(3)$.
		  \begin{itemize}
		   \item[$(1)$]$X_n = \supp \mu_n$ is connected.
		   \item[$(2)$]For any $r>0$ and $n\in \mathbb{N}$, all measures of closed balls in $X_n$ with radius $r$ are mutually equal.
		   \item[$(3)$]For each $n\in \mathbb{N}$, there exists an isometric embedding from $X_n$ to $Y_n$.
		  \end{itemize}
	 Then, for any $\kappa$ with $0<\kappa < \inf\limits_{n\in \mathbb{N}}m_n$, we have
	 \begin{align*}
	  \liminf_{n\to \infty} \diam (X_n \oblip Y_n,m_n -\kappa)\geq \frac{1}{2}\liminf_{n\to \infty} \diam X_n.
	  \end{align*}
	\end{prop}
	From Proposition \ref{kantan}, we obtain the following corollary:

	\begin{cor}For any $\kappa$ with $0<\kappa <1$, we have
	 \begin{align*}
	  \liminf_{n\to \infty} \diam (\mathbb{S}^n \oblip \mathbb{R}^{n+1},1-\kappa) >0.
	  \end{align*}
	 \end{cor}

	Note that the dimension $n+1$ of the screen can be replaced by any natural number greater than $n$.

	As an appendix, we discuss the case where the screen $Y$ is a (combinatorial) tree and give some answer to Exercise of
	Section $3\frac{1}{2}.32$ in \cite{gromov}. Precisely, we prove the following proposition:
		\begin{prop}\label{treenotoki}
	 Assume that a sequence $\{  (X_n,\dist_n,\mu_n ) \}_{n=1}^{\infty}$ of mm-spaces is a L\'{e}vy family. Then, we have
	 \begin{align*}
	  \sup \{\diam (X_n\oblip T,m_n-\kappa)    \mid T     \text{ is
	  a tree }               \}\to 0\ \text{as }n\to \infty
	  \end{align*}for any $\kappa>0$.
		 \end{prop}
	\section{Preliminaries}
		Let $(X,\dist)$ be a metric space. For $x\in X$, $r>0$, and $A,B \subseteq X$, we put
	\begin{align*}
	 B_X(x,r) :=\ &\{ y\in X \mid \dist (x,y)\leq r      \}, \
	 \dist(A,B)  := \inf \{ \dist(a,b) \mid a\in A, b\in B    \}, \\
	 \dist(x,A):= \ &\dist(\{x\},A), \
	 A_r :=\{ y\in X \mid   \dist(y,A)\leq r  \}, \ 
	 A_{+r} := \{ y\in X \mid \dist(y,A) < r   \}. 
	 \end{align*}
	We denote by 
	$m$ and $m_n$ the total measure of mm-spaces $X$ and $X_n$ respectively, and by $\supp \mu$ the support of a Borel measure $\mu$. 

	\subsection{Observable diameter and L\'{e}vy radius}

	In this subsection, we prove several results in \cite{gromov} because we find no proof anywhere. Let $(X,\dist ,\mu)$ be an mm-space and $f:X\to \mathbb{R}$ a Borel measurable function.
	 A number $a_0 \in
	 \mathbb{R}$  is called a \emph{pre-L\'{e}vy mean} of $f$ if 
	  $f_{\ast} (\mu) \big((-\infty , a_0]\big)\geq m/2 \text{ {\rm and }} 
	 f_{\ast}(\mu) \big([a_0,+ \infty)\big) \geq m/2$.
	We remark that $a_0$ does exist, but it is not unique for $f$ in general. Let $A_f\subseteq \mathbb{R}$ be the set of all pre-L\'{e}vy means of $f$. The proof of the following lemma is easy and we omit the
	proof.
	\begin{lem}[{cf.~\cite[Section $3\frac{1}{2}.19$]{gromov}}]\label{mean}$A_f$ is a closed bounded interval.
	 \end{lem}
	  The \emph{L\'{e}vy mean} $m_f$ of $f$ is defined by $m_f:=(a_f +b_f)/2$, where $a_f:=\min A_f$ and $b_f:=\max A_f$. For $\kappa >0$, we define the \emph{L\'{e}vy
	radius}, $\ler (X;- \kappa)$, as the infimum of $\rho >0$ such that every $1$-Lipschitz function $f:X\to \mathbb{R}$ satisfies
	$\mu (\{ x\in X \mid |f(x)-m_f|\geq \rho      \})\leq \kappa$.

		\begin{lem}[{cf.~\cite[Section $3\frac{1}{2}.32$]{gromov}}]\label{ccc}For any $\kappa >0$, we have
	 \begin{align*}
	   \diam(X  \oblip \mathbb{R},m-\kappa) \leq 2 \ler (X; -\kappa).
	  \end{align*}
	 \begin{proof}Put $\rho := \ler (X;-\kappa)$. It follows from the definition of the L\'{e}vy radius that $\mu (\{ x\in X
	  \mid |f(x)-m_f|   \geq \rho \}) \leq \kappa$ for any $1$-Lipschitz function $f:X\to \mathbb{R}$. So we obtain $f_{\ast}
	  (\mu) \big( (m_f-\rho, m_f + \rho ) \big)\geq m-\kappa$, which implies $\diam (f_{\ast}(\mu), m-\kappa)\leq \diam (m_f-\rho,m_f +\rho)
	  =2\rho$. This completes the proof.
	 \end{proof}
	\end{lem}
	\begin{lem}[{cf.~\cite[Section $3\frac{1}{2}.32$]{gromov}}]\label{bb}For any $\kappa$ with $0 < \kappa <m /2$ we have 
	 \begin{align*}
	 \ler (X; -\kappa) \leq
	 \diam(X\oblip \mathbb{R}, m-\kappa).
	 \end{align*}
	 \begin{proof}
	  Put $a:= \diam (X\oblip \mathbb{R},m-\kappa)$. For any $\varepsilon >0$, there exists a closed interval $X_0 \subseteq \mathbb{R}$
	  such that $f_{\ast}(\mu)(X_0)\geq m-\kappa$ and $ \diam X_0\leq a+\varepsilon$. We shall show that $m_f \in X_0$. 
	  If $X_0 \subseteq (-\infty ,m_f)$, we have
	  \begin{align*}
	   \frac{m}{2}<m-\kappa \leq f_{\ast}(\mu)(X_0)\leq f_{\ast}(\mu)\big((-\infty , m_f) \big)\leq \frac{m}{2},
	  \end{align*}
	  which is a contradiction. In the same way, we have $X_0 \nsubseteq (m_f,+\infty)$. 
	  Hence, we get $m_f \in  X_0$, which yields $X_0\subseteq [m_f-(a+\varepsilon),m_f + (a+\varepsilon)]$ since $\diam X_0 \leq
	  a+\varepsilon$. Therefore, we obtain
	  \begin{align*}
	   \mu (\{ x\in X \mid |f(x)-m_f| \leq a+\varepsilon \})=\ & f_{\ast}(\mu)([m_f-(a+\varepsilon),m_f+(a+\varepsilon)])\\
	   \geq \ & f_{\ast}(\mu)(X_0)\geq m-\kappa. 
	  \end{align*}As a result, we have $\ler (X; -\kappa)\leq a+\varepsilon$, which completes the proof of the lemma. 
	 \end{proof}
	\end{lem}

	Combining Lemma \ref{ccc} with Lemma \ref{bb}, we obtain the following corollary: 
	\begin{cor}[{cf.~\cite[Section $3\frac{1}{2}.32$]{gromov}}]A sequence $\{  X_n    \}_{n=1}^{\infty}$ of mm-spaces is a L\'{e}vy family if and only if 
	 $\ler (X_n; -\kappa)\to 0$ as $n\to \infty$ for any $\kappa >0$. 
	 \end{cor}

	\subsection{L\'{e}vy radius and concentration function}

	 Given an mm-space $(X,\dist ,\mu)$, we define the function $\alpha_{X}:(0,+\infty)\to \mathbb{R}$ by 
	  \begin{align*}
	   \alpha_{X}(r):= \sup \{\mu (X\setminus A_{+r}) \mid A \text{ {\rm is a Borel subset of }}X \text{ {\rm such that }} \mu (A) \geq m/2 \}, 
	  \end{align*}
	   and call it the \emph{concentration function} of $X$. Although the following lemmas and corollary are somewhat standard, we prove them for the completeness of this paper.
	\begin{lem}[{cf.~\cite[Section $1.3$]{ledoux}}]\label{a1}For any $r>0$ we have
	  \begin{align*}
	   \ler (X;-2\alpha_X(r))\leq r. 
	  \end{align*}
	  \begin{proof}Let $f:X\to \mathbb{R}$ be a $1$-Lipschitz function. We put $A:=\{ x\in X \mid f(x)\leq m_f \}$ and $A':=\{ x\in X \mid
	   m_f \leq f(x)\}$. Then,
	   \begin{align*}
		\{ x\in X \mid f(x)\geq m_{f}+ r \} = X \setminus \{ x\in X\mid f(x) <
		m_{f} + r\}
		\subseteq  X\setminus A_{+r}. 
	   \end{align*}
	   In the same way,
	   \begin{align*}
		\{ x\in X \mid m_f \geq f(x)+ r \} \subseteq X\setminus A'_{+r}. 
	   \end{align*}
	   Since $\mu(A)\geq m/2$ and $\mu(A')\geq m/2$, we have
	   \begin{align*}
		\mu (\{ x\in X \mid |f(x)-m_f |\geq r \}) 
		= \ &\mu(\{ x\in X \mid f(x) \geq m_f
		+ r\})+ \mu (\{ x\in X \mid m_f \geq f(x)+ r \})\\
		\leq \ &\mu(X\setminus  A_{+r})+
		\mu(X\setminus A'_{+r} )\\
		\leq \ &2 \alpha_X (r).
	   \end{align*}
	   This completes the proof. 
	  \end{proof}
	  \end{lem}

		  \begin{lem}[{cf.~\cite[Section 1.3]{ledoux}}]\label{a2}
	   For any $\kappa$ with $0<\kappa < m/2$, we have
	   \begin{align*}
		\alpha_X \big(2\ler(X;-\kappa)\big)\leq \kappa.
	   \end{align*}
	\begin{proof}Let $A$ be a Borel subset of $X$ such that $\mu(A)\geq m/2$. We define a 
		function $f:X\to \mathbb{R}$
		 by $f(x):=\dist(x,A)$. Putting $\rho:=\ler(X;-\kappa)$, by the definition of the L\'{e}vy radius, we have $\mu \big( \{ x\in X \mid |f(x)-m_f|\geq \rho  \} \big)\leq
		\kappa$. Then we get
		\begin{align*}
		 \mu\big(      \{x\in X \mid |f(x)-m_f|< \rho\} \cap A \big)\geq
		 \ & 
		 \mu \big(  \{x\in X \mid |f(x)-m_f|< \rho\} \big)+\mu(A)-m\\ 
		 \geq \ & (m-\kappa) +m/2-m=m/2 -\kappa>0. 
		 \end{align*}Hence, there exists a point $x_0 \in \{x\in X \mid |f(x)-m_f|<
		\rho\} \cap A$ and we have $m_f=|f(x_0)-m_f|<\rho$. Consequently, we obtain
		\begin{align*}
		 \mu(X\setminus A_{+ 2\rho})=\mu \big( \{ x\in X \mid f(x)\geq
		 2\rho   \}  \big)\leq \mu \big( \{ x\in X \mid   |f(x)-m_f|
		 \geq \rho  \} \big)\leq \kappa,
		 \end{align*}which completes the proof of the lemma.
	   \end{proof}
		   \end{lem}

	  \begin{cor}A sequence $\{ X_n \}_{n=1}^{\infty}$ of mm-spaces is a L\'{e}vy family if and only if 
	   $\alpha_{X_n}(r)\to 0$ as $n\to \infty$ for any $r>0$.
	   \begin{proof}Let $\{ X_n \}_{n=1}^{\infty}$ be a L\'{e}vy family. Fix $r>0$ and take any $\varepsilon >0$. For an $n\in \mathbb{N}$ with $m_n/2 \leq
		\varepsilon$, we have $\alpha_{X_n}(r)\leq
		\varepsilon$. Hence, we only consider the case of $m_n/2 >\varepsilon$. From the assumption, we have
		$2\ler(X_n;-\varepsilon) \leq r$ for any sufficiently
	    large $n\in \mathbb{N}$. Therefore, by virtue of Lemma \ref{a2}, we have
		\begin{align*}
		 \alpha_{X_n}(r)\leq \alpha_{X_n}\big( 2\ler(X_n;-\varepsilon)   \big)\leq \varepsilon, 
		 \end{align*}which shows $\alpha_{X_n}(r)\to 0$ as $n\to \infty$.

		Conversely, assume that $\alpha_{X_n}(r)\to 0$ as $n\to \infty$ for any $r>0$. Fix $\kappa > 0$ and take any $\varepsilon >
		0$. From the assumption, we have $2\alpha_{X_n}(\varepsilon)\leq \kappa$ for any sufficiently large $n\in \mathbb{N}$. Therefore, applying Lemma
		\ref{a1} to $X_n$, we obtain
		\begin{align*}
		 \ler (X_n;-\kappa)\leq \ler \big(X_n;-2\alpha_{X_n}(\varepsilon)\big)\leq \varepsilon.
		 \end{align*}This completes the proof. 
		\end{proof}
	 \end{cor}

		 \subsection{Concentration function and separation distance}
	Let $(X, \dist , \mu )$ be an mm-space. For any $\kappa_0, \kappa_1, \cdots
	 ,\kappa_N \in \mathbb{R}$, we define 
	 \begin{align*}
	  \sep (X ; \kappa_0 , \cdots ,\kappa_N)=\ & \sep (\mu ; \kappa_0, \cdots
	  , \kappa_N)\\
	 : =\ & \sup \{  \min_{i\neq j} \dist(X_i, X_j)   \mid   X_0 , \cdots , X_N
	   \text{ {\rm are Borel subsets of }} X \\ & \hspace{4.9cm}\text{{\rm which satisfy }}
	  \mu(X_i) \geq \kappa_i \text{ {\rm for any }}i \}, 
	 \end{align*}and call it the \emph{separation distance} of $X$. In this subsection, we investigate relationships between
	the concentration function and the separation distance. The proof of the following lemma is easy, and we omit the proof.
	\begin{lem}[{cf.~\cite[Section $3\frac{1}{2}.33$]{gromov}}]\label{neko}
	 Let $(X, \dist_X, \mu_X)$ and $(Y, \dist_Y , \mu_Y) $ be two mm-spaces. Assume that a $1$-Lipschitz map $f : X\to Y$
	 satisfies $f_{\ast}(\mu_X)= \mu_Y$. Then we have
	 \begin{align*}
	  \sep (Y; \kappa_0, \cdots , \kappa_N) \leq \sep (X ; \kappa_0, \cdots , \kappa_N).
	 \end{align*}
	\end{lem}
	Let us recall that the \emph{Hausdorff distance} between two bounded closed subsets $A$ and $B$ in a metric space $X$ is
	defined by
	\begin{align*}
	 \dist_H (A,B):= \inf \{ \varepsilon > 0 \mid  A \subseteq B_{+\varepsilon}, \ B \subseteq A_{+\varepsilon}    \}. 
	 \end{align*}
	It is easy to check that $\dist_H$ is the metric on the set $\mathcal{C}_X$ of all bounded closed subsets of $X$. 
	\begin{lem}[{Blaschke, cf.~\cite[Theorem 4.4.15]{ambro}}]\label{compact}
	 If $X$ is a compact metric space, then $(\mathcal{C}_X, \dist_H)$ is also compact. 
	\end{lem}

	\begin{lem}Let $(X,\dist,\mu)$ be an mm-space and assume that $\supp \mu$ is connected. Then, for any $r>0$ with
	 $\alpha_X(r)>0$ 
	  we have 
	 \begin{align*}
	  \sep \Big(X;\frac{m}{2},\alpha_X(r)\Big)\leq r.
	  \end{align*}
	 \begin{proof}The proof is  by contradiction. We may assume that $X=\supp \mu$. Suppose that $\sep \big(X;m/2,\alpha_X(r)\big)> r$,
	  there exist 
	  $r_0 >0$ with $r_0 >r$ and Borel subsets $X_1, X_2 \subseteq X$ such that
	  $\mu(X_1)\geq m/2$, $\mu(X_2)\geq \alpha_X(r)$, and $\dist(X_1,X_2)>r_0$. Let us show that $(X_1)_{r} \nsupseteq (X_1)_{+
	  r_0}$. If $(X_1)_r = (X_1)_{+r_0}$, we have 
	  $X=(X_1)_r \cup \big(X \setminus (X_1)_{+ r_0}  \big)$. Since $X$ is connected, we get either $X_1= \emptyset$ or $X
	  \setminus (X_1)_{+ r_0}   =\emptyset$. It follows from $\mu
	  \big((X_1)_r\big)\geq \mu(X_1)\geq m/2 >0$ that $(X_1)_r \neq  \emptyset$
	  . By $\dist(X_1,X_2)>r_0$, we obtain $X_2 \subseteq X\setminus (X_1)_{+ r_0}$, which implies that $\mu \big( X\setminus (X_1)_{+ r_0}
	  \big)\geq \mu(X_2) \geq \alpha_X(r)>0$. Therefore, we have $X\setminus (X_1)_{+ r_0} \neq \emptyset$, which is a
	  contradiction. Thus, there exists a point $x_0 \in \big(X\setminus (X_1)_r \big) \setminus
	  \big( X\setminus (X_1)_{+ r_0}  \big)$. Taking a sufficiently small ball $B$ centered at $x_0$ such that $B\subseteq X\setminus
	  (X_1)_r$ and $B\cap \big( X\setminus (X_1)_{+ r_0}\big) =\emptyset$, we have 
	  \begin{align*}
	   \mu(X_2)\geq \alpha_X(r)\geq \ &\mu   \big(X\setminus (X_1)_r \big) \\
	   \geq \ & \mu    \big(  B \cup (X\setminus (X_1)_{+ r_0})
	   \big)=\mu(B)+\mu \big( X\setminus (X_1)_{+ r_0}\big)\\
	   > \ & \mu \big( X\setminus (X_1)_{+ r_0}       \big)\geq \mu(X_2),
	   \end{align*}which is a contradiction. Therefore, we have
	  finished the proof.
	  \end{proof}
	 \end{lem}

	\begin{rem}
	 If $\ \supp \mu$ is disconnected, 
	 the above lemma does not hold in general. For example, consider the space $X:=\{ x_1 , x_2\}$ with a metric $\dist$ given by
	 $\dist (x_1,x_2):=1$ and with a Borel probability measure $\mu$ given by $\mu (\{ x_1\})=\mu (\{  x_2   \}):= 1/2$. In this
	 case, we have $ \alpha_{X}(1/2)=1/2$ and $ \sep (\mu ,1/2 ,1/2)=1 $. 
	 \end{rem}

	\begin{lem}\label{a9}
	 For any $r>0$ there exists a Borel subset $X_0\subseteq X$ such that 
	 \begin{align*}
	  \mu \big( X\setminus (X_0)_{+r}\big)= \alpha_X(r) \text{ and } \mu(X_0)\geq \frac{m}{2}.
	  \end{align*}
	 \begin{proof}From the definition of the concentration function, for any $n\in \mathbb{N}$, there exist a closed subset $A_n
	  \subseteq X$ such that
	  \begin{align*}
	   \mu(A_n)\geq \frac{m}{2} \text{ and } \mu \big( X \setminus (A_n)_{+r}\big) +\frac{1}{n} \geq \alpha_{X}(r).
	   \end{align*}Take an increasing sequence $K_1 \subseteq K_2 \subseteq \cdots $ of compact subsets of $X$ such that
	  $\mu(K_n)\to m$ as $n \to \infty$. By using Lemma \ref{compact} and the diagonal argument, we have that $\{
	  A_n\cap K_i\}_{n=1}^{\infty}$ Hausdorff converges to a closed subset $B_i \subseteq K_i$ and $\{
	  \big(X\setminus (A_n)_{+r} \big) \cap K_i \}_{n=1}^{\infty}$ Hausdorff converges to a closed subset $C_i
	  \subseteq K_i$ for each $i\in \mathbb{N}$. Put $K:=\bigcup_{i=1}^{\infty}A_i, \widetilde{K}:=
	  \bigcup_{i=1}^{\infty}C_i$, $X_0:= \overline{K}$, and $Y_0:= \overline{\widetilde{K}}$. It is easy to check that $B_1\subseteq
	  B_2 \subseteq \cdots$ and $C_1\subseteq C_2
	  \subseteq \cdots$. We will show that $\dist (K,\widetilde{K})\geq r$ by contradiction. If $\dist(K,\widetilde{K})<r$, there exists $r_0>0$
	  such that $\dist(K,\widetilde{K})<r_0 < r$. Hence there exist $x\in K$ and $y\in \widetilde{K}$ such that
	  $\dist(x,y)<r_0$. There exists $i\in \mathbb{N}$ such that $x\in B_i$ and $y\in C_i$, because both $\{B_n\}_{n=1}^{\infty}$ and $\{C_n
	  \}_{n=1}^{\infty}$ are increasing sequences. Since both $\{
	  A_n\cap K_i   \}_{i=1}^{\infty}$ and $\{  \big(X\setminus  (A_n)_{+r}\big)\cap K_i     \}_{n=1}^{\infty}$
	  Hausdorff converge to $B_i$ and $C_i$ respectively, 
	  there exist two sequences $\{ x_n  \}_{n=1}^{\infty}$, $\{ y_n   \}_{n=1}^{\infty} \subseteq X$ such that $\dist(x_n,x),
	  \dist(y_n,y)\to 0$ as $n\to \infty$, 
	  and $x_n\in A_n, \dist(y_n,A_n)\geq r$ for any $n\in \mathbb{N}$. Therefore, for any sufficiently large $n\in \mathbb{N}$ we have
	  \begin{align*}
	   \dist(x_n,y_n)\leq \dist(x_n,x)+\dist(x,y)+\dist(y_n,y)<r_0+\dist(x_n,x)+\dist(y_n,y)<r,
	   \end{align*}which is a contradiction, because $x_n \in A_n$ and $\dist(y_n,A_n)\geq r$. Thus, we obtain $\dist(K,\widetilde{K})\geq r$
	  which yields that $\dist (X_0,Y_0)\geq r$ and therefore $Y_0 \subseteq X\setminus (X_0)_{+r}$. Let us show that
	  $\mu(X_0)\geq m/2$ and $
	  \mu(Y_0)\geq \alpha_{X}(r)$. For any $\varepsilon >0$ there exists 
	  $n_0\in \mathbb{N}$ such that $\mu(K_{n_0})+\varepsilon \geq m$. Take any $\delta >0$. Then, for any sufficiently large $m$
	  we have $A_m \cap K_{n_0} \subseteq (B_{n_0})_{\delta}$. Therefore,
	  \begin{align*}
	   \mu\big((X_0)_{\delta} \big)\geq \mu \big( (B_{n_0})_{\delta} \big)\geq \mu(A_m \cap K_{n_0})\geq
	   \mu(A_m)+\mu(K_{n_0})-m\geq \frac{m}{2}-\varepsilon.
	   \end{align*}By taking $\delta \to 0$, we obtain $\mu(X_0)\geq m/2-\varepsilon$, which shows $\mu(X_0)\geq m/2$. In the same
	  way, we have $\mu (Y_0)\geq \alpha_X(r)$. This completes the proof. 
	  \end{proof}
	 \end{lem}

	 Lemma \ref{a9} directly implies 
	\begin{lem}\label{a10}For an mm-space $X$ and $r>0$, we have 
	 \begin{align*}
	  \sep \Big(X;\frac{m}{2},\alpha_X(r)\Big)\geq r.
	  \end{align*}
	 \end{lem}

	 \begin{cor}\label{dd}
	  If a sequence $\{X_n \}_{n=1}^{\infty}$ of mm-spaces satisfies that $\sep(X_n;\kappa,\kappa)\to 0$ as $n \to
	  \infty$ for any $\kappa >0$, we then have $\alpha_{X_n}(r)\to 0$ as $n\to \infty$ for any $r>0$. 
	  \begin{proof}
	   Suppose that there exists $c>0$ such that 
	   $\alpha_{X_n}(r)\geq c$ for infinitely many $n\in \mathbb{N}$. Applying Lemma 
	   \ref{a10} to $X_n$, we have
	   \begin{align*}
		r\leq \sep\Big(X_n;\frac{m_n}{2},\alpha_{X_n}(r)\Big) \leq \sep
		\big(X_n ; \alpha_{X_n}(r),\alpha_{X_n}(r)  \big) \leq \sep (X_n;c,c).
		\end{align*}This is a contradiction, since the right-hand side of the above inequality converges to $0$ as $n \to \infty$. This
	   completes the proof. 
	  \end{proof}
	 \end{cor}

	 \begin{lem}[{cf.~\cite[Lemma $1.1$]{ledoux}}]\label{cc}Let
	  $(X,\dist,\mu)$ be an mm-space. Assume that a Borel subset $A\subseteq X$ and 
	  $r_0 >0$ satisfy $\mu(A)\geq \kappa$ and $\alpha_X(r_0)<\kappa$. Then, for any $r>0$ we have 
	  \begin{align*}
	   \mu\big( X\setminus A_{+(r_0+r)} \big)\leq \alpha_X(r).
	   \end{align*}
	  \end{lem}

	  \begin{cor}\label{ee}Assume that a sequence $\{  (X_n,\dist_n,\mu_n)    \}_{n=1}^{\infty}$ of mm-spaces satisfy
	   $\alpha_{X_n}(r) \to 0$ as $n\to \infty$ for any $r>0$. Then, we have $\sep (X_n;\kappa,\kappa)\to 0$ as $n\to \infty$ for
	   any $\kappa >0$.
	   \begin{proof}Since $\sep(X_n;\kappa,\kappa)=0$ for $n\in \mathbb{N}$ with $m_n< \kappa$, we assume that $m_n\geq \kappa$ for any
		$n\in \mathbb{N}$. For any $\varepsilon >0$, we have $\alpha_{X_n}(\varepsilon)< \kappa /2$ for any sufficiently
		large $n\in \mathbb{N}$ from the assumption. Thus, it follows from Lemma \ref{cc} that $\mu_{n}\big(  X_n \setminus (A_n)_{+2\varepsilon}    \big)\leq
		\alpha_{X_n}(\varepsilon)< \kappa/2$ for any Borel sets $A_n,B_n\subseteq X_n$ with $\mu_{n}(A_n), \mu_{n}(B_n)\geq
		\kappa$. In the same way, we get $\mu_{n}\big(  X_n\setminus (B_n)_{+2\varepsilon}   \big)<\kappa /2$. Therefore, we obtain
		\begin{align*}
		 \mu_{n}\big(X_n \setminus \big(  (A_n)_{+2\varepsilon}   \cap (B_n)_{+2\varepsilon}     \big)       \big) \leq \ &
		 \mu_{n}\big(  X_n \setminus (A_n)_{+2\varepsilon} \big) + \mu_{n}\big(  X_n \setminus (B_n)_{+2\varepsilon}
		 \big)\\
		 < \ & \kappa \leq m_n,
		 \end{align*}which implies $\mu_{n}\big((A_n)_{+2\varepsilon}\cap (B_n)_{+2\varepsilon}\big) >0$, and thereby
		$(A_n)_{+2\varepsilon}\cap (B_n)_{+2\varepsilon} \neq \emptyset$. As a consequence, we have $\diam(A_n,B_n)\leq
		4\varepsilon$, which shows 
		$\sep (X_n;\kappa,\kappa)\leq 4\varepsilon$. This completes the proof. 
		\end{proof}
	   \end{cor}

	Combining Corollary \ref{dd} with Corollary \ref{ee}, we obtain the following corollary:
	\begin{cor}[{cf.~\cite[Section $3\frac{1}{2}.33$]{gromov}}]A sequence $\{ X_n \}_{n=1}^{\infty}$ of mm-spaces is a L\'{e}vy family if and only if 
	 $\sep(X_n;\kappa,\kappa)\to 0$ as $n\to \infty$ for any $\kappa >0$.
	 \end{cor}
 \subsection{$L^p$-concentration of maps} 
In this subsection, we investigate relationships between concentration
of maps and $L^p$-concentration of maps. The standard argument yields the following basic fact.
\begin{lem}\label{kouta1}For $0<p<q$, $L^q$-concentration of maps implies
 $L^p$-concentration of maps under the assumption $\sup\limits_{n\in
 \mathbb{N}} m_n < +\infty$.
\end{lem}

\begin{lem}\label{var1}Let $X$ be an mm-space and $Y$ a metric space. For any $\kappa>0$, $p\in (0,+\infty]$, and a Borel measurable map $f:X\to
 Y$, we have
\begin{align*}
\diam (f_{\ast}(\mu_X),m-\kappa)\leq \frac{2}{(\kappa m)^{1/p}}V_p(f).
\end{align*}
\begin{proof}The case of $p=\infty$ is easy, so we consider the case of
 $p<+\infty$ only.
Assume that $\mu_X \big(\big\{ x\in X \mid \dist_Y\big( f(x),f(x') \big)
 \geq \varepsilon_0                              \big\}\big)>\kappa$
 for any $x'\in X$ and $\varepsilon_0:= V_p(f)/(\kappa m)^{1/p}$. By the Chebyshev's inequality, we get
\begin{align*}
V_p(f)^p= \int_X \Big\{ \int_{X}
 \dist_{Y}\big(f(x),f(x')\big)^p       \ d\mu_X(x)       \Big\} \
 d\mu_X(x')
>  \int_X \varepsilon_0^p \kappa \ d \mu_X(x') = \varepsilon_0^p \kappa
 m = V_p(f)^p.
\end{align*}Hence, there is a point $x'\in
 X$ such that $\mu_X\big(\big\{ x\in X \mid \dist_Y\big( f(x),f(x') \big)
 \geq \varepsilon_0                              \big\}\big)\leq
 \kappa$. This completes the proof.
\end{proof}
\end{lem}

\begin{cor}\label{hayakuowaritai1}For any $p\in (0, +\infty]$, $L^p$-concentration of maps
 implies concentration of maps. 
\end{cor}

\begin{lem}\label{var2}Let $X$ be an mm-space with a finite diameter and $Y$ a
 metric space. Then, for any $\kappa>0$, $p\in (0,+\infty)$, and
 $1$-Lipschitz map $f:X\to Y$, we have
\begin{align*}
V_p(f)^p \leq m^2 \diam
 \Big(f_{\ast}(\mu_X),m-\frac{\kappa}{(\diam X)^p}\Big)^p + \Big(
 2m-\frac{\kappa}{(\diam X)^p}\Big)\kappa.
\end{align*}
\begin{proof}For any $\varepsilon >0$ with $\diam
 \big(f_{\ast}(\mu_X),m-\kappa/(\diam X)^p\big) <\varepsilon$, there
 exists a Borel subset $A\subseteq Y$ such that $\diam A <\varepsilon$
 and $f_{\ast}(\mu_X)(A)\geq m-\kappa / (\diam X)^p$. By $\diam A <\varepsilon$, we get
\begin{align}\label{miyukichi1}
\int \int_{f^{-1}(A)\times f^{-1}(A)} \dist_Y\big( f(x),f(x') \big)^p \
 d\mu_X(x) d\mu_X(x') \leq m^2 \varepsilon^p. 
\end{align}Since
\begin{align*}
(\mu_X \times \mu_X)\big(X\times X \setminus f^{-1}(A)\times
 f^{-1}(A)\big)\leq \ & m^2- \Big(   m- \frac{\kappa}{(\diam X)^p}
 \Big)^2\\ = \ & \Big(
 2m-\frac{\kappa}{(\diam X)^p}\Big)\frac{\kappa}{(\diam X)^p},
\end{align*}we also have
\begin{align}\label{miyukichi}
\int \int_{X\times X \setminus f^{-1}(A)\times f^{-1}(A)} \dist_Y
 \big(f(x),f(x')\big)^p \ d\mu_X(x) d\mu_X(x') \leq \Big(
 2m-\frac{\kappa}{(\diam X)^p}\Big)\kappa. 
\end{align}Combining (\ref{miyukichi1}) with (\ref{miyukichi}), we obtain
\begin{align*}
V_p(f)^p \leq  m^2 \varepsilon^p +\Big(
 2m-\frac{\kappa}{(\diam X)^p}\Big)\kappa. 
\end{align*}This completes the proof.
\end{proof}
\end{lem}

\begin{cor}\label{waruiyokan}Let $\{ X_n     \}_{n=1}^{\infty}$ be a sequence of mm-spaces
 with finite diameters and $\{  Y_n\}_{n=1}^{\infty}$ a sequence of
 metric spaces. Assume that $\sup\limits_{n\in \mathbb{N}}m_n < +\infty$
 and 
\begin{align*}
\diam \Big(X_{n} \oblip Y_n, m_{n}- \frac{\kappa}{  (\diam X_{n})^p}\Big)
		 \to 0 \text{ as }n\to \infty
\end{align*}for any $\kappa >0$ and some $p\in (0,+\infty)$. Then we
 have $\obin_{Y_n}(X_n) \to 0$ as $n\to \infty$.
\end{cor}

\subsection{Basics of the barycenter of a measure}
First, we shall review some standard terminologies in metric geometry. A metric
space $(X,\dist_X)$ is said to be \emph{geodesic} if every two points
$x,y \in X$ can be joined by a rectifiable curve whose length coincides
with $\dist_X(x,y)$. A rectifiable curve $\gamma:[0,1]\to X$ is called a
\emph{minimal geodesic} if it is globally minimizing and has a constant speed,
i.e., parameterized proportionally to the arclength. A geodesic metric
space $X$ is called a CAT$(0)$-space if we have
\begin{align*}
\dist_X\big(x,\gamma (1/2)\big)^2 \leq \frac{1}{2}\dist_X(x,y)^2 +
 \frac{1}{2}\dist_X(x,z)^2 - \frac{1}{4} \dist_X(y,z)^2
\end{align*}for any $x,y,z \in X$ and any minimizing geodesic $\gamma
 :[0,1]\to X$ from $y$ to $z$. If the reverse inequality
\begin{align*}
\dist_X\big(x,\gamma (1/2)\big)^2 \leq \frac{1}{2}\dist_X(x,y)^2 +
 \frac{1}{2}\dist_X(x,z)^2 - \frac{1}{4} \dist_X(y,z)^2
\end{align*}holds for any $x,y,z$ and $\gamma$, then we say that
 $(X,\dist_X)$ is an \emph{Alexandrov space of nonnegative curvature}.
\begin{ex}\upshape 
$(1)$ A complete simply connected Riemannian manifold is
	    a CAT(0)-space if and only if its sectional curvature is
	    nonpositive everywhere. Similarly, a complete Riemannian
	    manifold is an Alexandrov space of nonnegative curvature if
	    and only if its sectional curvature is nonnegative.

$(2)$ A Hilbert spaces is both CAT(0)-space and Alexandrov space of
	    nonnegative curvature.

$(3)$ For a compact convex domain $\Omega \subseteq \mathbb{R}^n$,
	    let $ \partial \Omega$ equip the length metric
	    $\dist$
	    induced from the standard metric of $\mathbb{R}^n$. Then
	    $(\partial \Omega ,\dist)$ is an Alexandrov space of
	    nonnegative curvature.

$(4)$ A tree is a CAT(0)-space. 
\end{ex}
See \cite{bri}, \cite{bur}, \cite{bur2}, \cite{gromov}, \cite{jost} and references
therein for basics of CAT(0) spaces and Alexandrov spaces of nonnegative curvature.

Let $(X,\dist_X)$ be a metric space. We denote by $B(X)$ the set of all
finite Borel measures on $X$ with
separable supports. Given $1\leq p < +\infty$, we indicate by
$B^p(X)$ the set of all Borel measures $\nu \in B(X)$ such that $\int_{X}\dist_X(x,y)^p \ d\nu (y)<+\infty$ for
some (hence all) $x\in X$. Obviously, $B^p(X)\subseteq B^q(X) \subseteq B(X)$ for any $p>q$.

For any $\nu \in B^1(X)$ and $z\in X$, we consider the function
$h_{z,\nu}:X\to \mathbb{R}$ defined by
\begin{align*}
h_{z,\nu}(x):= \int_{X} \{\dist_X(x,y)^2-\dist_X(z,y)^2 \} \ d\nu(y).
\end{align*}
Note that
\begin{align*}
\int_X |\dist_X(x,y)^2-\dist_X(z,y)^2| \ d\nu(y)\leq \dist_X(x,z)\int_X
 \{\dist_X(x,y)+ \dist_X(z,y)\} \ d \nu(y) <+\infty,
\end{align*}that is, the function $h_{z,\nu}$ is well-defined.

\begin{dfn}\upshape The point $z_0 \in X$ is called the
 \emph{barycenter} of $\nu \in B^1(X)$ if for any $z\in X$, $z_0$ is the unique minimizing
 point of the function $h_{z,\nu}$. We denote the point $z_0$ by $b(\nu)$.
A metric space $X$ is said to be \emph{barycentric} if every $\nu \in B^1(X)$ has the barycenter. 
\end{dfn}

\begin{prop}[{cf.~\cite[Proposition
 $4.3$]{sturm}}]\label{onegai}
A CAT(0)-space $N$ is barycentric. Moreover, for any $\nu \in B^2(N)$,
 we have
\begin{align*}
\int_N \dist_N(b(\nu),y)^2 \ d\nu(y) = \inf_{x\in N} \int_{N}
 \dist_N(x,y)^2 \ d\nu(y).
\end{align*}
\end{prop}
S. Ohta pointed out in \cite[Lemma $5.5$]{ohta} that Proposition
\ref{onegai} holds for more general metric spaces, such as $2$-uniformly
convex metric spaces.

A simple variational arguments yields the following two lemmas.
\begin{lem}[{cf.~\cite[Propsition $5.4$]{sturm}}]\label{atata}
Let $H$ be a Hilbert space. Then for each $\nu \in B^1(H)$ with $m=\nu(X)$, we have
\begin{align*}
b(\nu)= \frac{1}{m}\int_H y \ d\nu(y).
\end{align*}
\end{lem}

\begin{lem}[{cf.~\cite[Proposition $5.10$]{sturm}}]\label{atatata}Let $N \in
 \mathcal{NM}^n$ and $\nu \in B^1(N)$. Then $x=b(\nu)$ if and only if
\begin{align*}
\int_N \exp_x^{-1} (y) \ d\nu (y)=0.
\end{align*}In particular, $b\big((\exp_{b(\nu)}^{-1})_{\ast}(\nu)\big)=0.$
\end{lem}

\subsection{Some property of observable diameter $\diam (X\oblip \mathbb{R}^k,m-\kappa)$}
For a measure space $(X,\mu)$ with $\mu(X)<+\infty$, we denote by
$\mathcal{F}(X,\mathbb{R}^k)$ the space of all maps from $X$ to
$\mathbb{R}^k$. Given $\lambda\geq 0$ and $f,g\in \mathcal{F}(X,\mathbb{R}^k)$, we put 
\begin{align*}
\me_{\lambda}(f,g):= \inf \{  \varepsilon >0 \mid \mu 
 \big(\{ x\in X \mid |f(x)-g(x)|\geq \varepsilon \} \big)\leq \lambda \varepsilon             \}.
\end{align*}Note that this $\me_{\lambda}$ is a distance function on $\mathcal{F}(X,
 \mathbb{R}^k)$ for any $\lambda \geq 0$ and its topology on
 $\mathcal{F}(X, \mathbb{R}^k)$ coincides with the topology of the 
convergence in measure for any $\lambda >0$. Also, the distance functions
 $\me_{\lambda}$ for all $\lambda >0$ are mutually equivalent. 

Let $\lambda \geq 0$ and $\varepsilon>0$. A map from an mm-space to a metric space, say
$f:X\to Y$ is called \emph{$\lambda$-Lipschitz up to $\varepsilon$} if 
\begin{align*}
\dist_Y \big(f(x),f(x')\big) \leq \lambda \dist_X (x,x') +\varepsilon
\end{align*}for all $x,x'$ in a Borel subset $X_0\subseteq X$ with
$\mu_X(X\setminus X_0)\leq \varepsilon$. 

\begin{prop}[{cf.~\cite[Section $3\frac{1}{2}.15$, $(3_b)$]{gromov}}]\label{kumogaippai}Let $(X,\dist_X, \mu_X)$, $(Y,\dist_Y, \mu_Y)$ be mm-spaces
 and $\lambda\geq 0$. Let $\varepsilon_n>0$ and $f_n:X \to Y$
 a $\lambda$-Lipschitz up to $\varepsilon_n$ Borel merasurable map and
 assume that $\varepsilon_n \to 0$ as $n\to \infty$ and the sequence $\{
 (f_n)_{\ast}(\mu_X) \}_{n=1}^{\infty}$ converges weakly to
 $\mu_Y$. Then, the sequence $\{ f_n \}_{n=1}^{\infty}$ has a
 $\me_1$-convergent subsequence.
\end{prop}
See \cite[Proposition $3.1$]{funafuna} for the details of proof of Proposition \ref{kumogaippai}.

For an mm-space $X$, we fix a point $x_0 \in \supp \mu_X$. Although the
following lemma is stated in \cite{gromov}, we give the proof of the lemma for the
completeness of this paper.
\begin{lem}[{cf.~\cite[Section $3.\frac{1}{2}.36$]{gromov}}]\label{mukatukubaka}Let $\{ f_n:B_X(x_0,n) \to
 \mathbb{R}^k\}_{n=1}^{\infty}$ a sequence of $1$-Lipschitz maps. Then,
 there exist a $1$-Lipschitz map $f:X\to \mathbb{R}^k$ and subsequence
 $\{f_{a(n)}\}_{n=1}^{\infty}\subseteq \{f_n\}_{n=1}^{\infty}$ such that
\begin{align*}
\me_1\big(f_{a(n)}-f_{a(n)}(x_0),f|_{B_X(x_0,a(n))}\big)\to 0 \text{ as }n\to \infty.
\end{align*}
\begin{proof}We only prove the case of $k=1$. The general case follows
 from this special case. Put
 $\mu_n:=\big(f_{n}-f_{n}(x_0)\big)_{\ast}(\mu_X|_{B_X(x_0,n)})$. We
 shall show that there exist a monotone increasing sequence 
 $\{a(n)\}_{n=1}^{\infty}\subseteq \mathbb{N}$ and a Borel measure $\nu$ on
 $\mathbb{R}$ such that the sequence
 $\{\mu_{a(n)}\}_{n=1}^{\infty}$
 weakly converges to $\nu$. Combining this with Proposition
 \ref{kumogaippai}, we get the proof of the lemma. Observe that
 $\mu_n([-m,m])\geq \mu_X\big(B_X(x_0,m)\big)>0$ for any $m,n\in
 \mathbb{N}$ with $n\geq m$. Since $[-1,1]$ is a
 compact metric space, by virtue of classical Prohorov's theorem
 (cf.~\cite[Theorem 5.1]{bil}), there exist a monotone increasing sequence $\{
 {a_1(n)}\}_{n=1}^{\infty} \subseteq \mathbb{N}$ and a
 Borel measure $\nu_1$ on $[-1,1]$ such that the sequence $\{
 \mu_{a_1(n)}|_{[-1,1]} \}_{n=1}^{\infty}$ converges weakly to
 the measure $\nu_1$. In the same way, there exist a monotone increasing
 subsequence $\{
 a_2(n)\}_{n=1}^{\infty} \subseteq \{  a_1(n)\}_{n=1}^{\infty}$ and a
 Borel measure $\nu_2$ on $[-2,2]$ such
 that the sequence $\{  \mu_{a_2(n)} \}_{n=1}^{\infty}$ converges weakly
 to the measure $\nu_2$. Observe that $\nu_2|_{[-1,1]}=\nu_1$. In this
 way, we obtain a monotone increasing sequence $\{
 a_{k}(n)\}_{n=1}^{\infty}\subseteq \{  a_{k-1}(n)\}_{n=1}^{\infty}$ and
 a Borel measure $\nu_m$ on $[-m,m]$. Define a Borel measure $\nu$ on
 $\mathbb{R}$ by $\nu|_{[-m,m]}:=\nu_{m}$ and put $a(n):=a_n(n)$. 

Take any open subset $O\subseteq \mathbb{R}$. For any $\varepsilon >0$,
 there is $m\in \mathbb{N}$ such that $\nu\big((-m,m)\cap
 O\big)+\varepsilon >\nu(O)$. Since $\nu([-m,m])/\mu_{a(n)}([-m,m])\to
 1$ as $n\to \infty$ and the sequence $\{\mu_{a(n)}|_{[-m,m]}\}_{n=1}^{\infty}$ converges
 weakly to the measure $\nu|_{[-m,m]}$, we get
\begin{align*}
\liminf_{n\to \infty}\mu_{a(n)}(O)\geq \liminf_{n\to
 \infty}\mu_{a(n)}\big((-m,m)\cap O\big) \geq \nu\big((-m,m)\cap O\big)>\nu(O)-\varepsilon.
\end{align*}Hence the sequence $\{\mu_{a(n)}\}_{n=1}^{\infty}$ converges
 weakly to the measure $\nu$. This completes the proof.
\end{proof}
\end{lem}
Put $m_n:=\mu_X\big(B_X(x_0,n)\big)$. 
\begin{cor}\label{kawaisouni2}For any $\kappa'>\kappa>0$, we have
\begin{align*}
\liminf_{n\to \infty}\diam (B_X(x_0,n)\oblip \mathbb{R}^k,m_n-\kappa')
 \leq \diam (X\oblip \mathbb{R}^k,m-\kappa).
\end{align*}
\begin{proof}Suppose that 
\begin{align*}
\liminf_{n\to \infty}\diam (B_X(x_0,n)\oblip \mathbb{R}^k,m_n-\kappa')
 >\alpha > \diam (X\oblip \mathbb{R}^k,m-\kappa).
\end{align*}Then, there is a sequence $\{ f_n:B_X(x_0,n)\to
 \mathbb{R}^k\}_{n=1}^{\infty}$ of $1$-Lipschitz maps such that $\diam
 ((f_n)_{\ast}(\mu_X),m_n-\kappa')> \alpha$ for any sufficiently large $n\in \mathbb{N}$. According to
 Lemma \ref{mukatukubaka}, there exist a subsequence
 $\{f_{a(n)}\}_{n=1}^{\infty}\subseteq \{ f_n\}_{n=1}^{\infty}$ and a
 $1$-Lipschitz map $f:X\to \mathbb{R}^k$ such that 
\begin{align}\label{kawaisouni}
\me_1\big(f_{a(n)}-f_{a(n)}(x_0),f|_{B_X(x_0,a(n))}\big)\to 0 \text{ as }n\to \infty.
\end{align}Assume that a closed subset $A\subseteq \mathbb{R}^k$
 satisfies that $f_{\ast}(\mu_X)(A)\geq m-\kappa$. Take $\varepsilon >0$
 with $2\varepsilon < \alpha -\diam (X\oblip \mathbb{R}^k,m-\kappa)$ and put 
\begin{align*}
A_n:= \big\{  x\in B_X\big(x_0,a(n)\big)   \mid |f_{a(n)}(x)-f_{a(n)}(x_0)-f(x)|<\varepsilon   \big\}.
\end{align*}Then, from (\ref{kawaisouni}), we get
\begin{align*}
(f_{a(n)})_{\ast}(\mu_X) \big((A)_{\varepsilon}+f_{a(n)}(x_0)\big)\geq
 f_{\ast}(\mu_X)(A)-\mu_X(X\setminus A_n)\geq m_n-\kappa'
\end{align*}for any sufficiently large $n\in \mathbb{N}$. Hence, we
 get $\diam \big((f_{a(n)})_{\ast}(\mu_X),m_n-\kappa'\big)\leq \diam A
 +2\varepsilon$ for any sufficiently large $n\in \mathbb{N}$, which
 implies that
\begin{align*}
\alpha< \diam \big((f_{a(n)})_{\ast}(\mu_X),m_n-\kappa'\big)\leq \diam
 (f_{\ast}(\mu_X),m-\kappa)+2\varepsilon  <\alpha
\end{align*}for any sufficiently large $n\in \mathbb{N}$. This is a
 contradiction and hence, the proof is completed.
\end{proof}
\end{cor}

	\section{Main theorems}

	Let $M$ be a complete Riemannian manifold and $\nu$ a Borel
	measure on $M$ such that $\supp \nu$ is compact and $\nu
	(M)<+\infty$. For $t>0$, we consider the function $\dist_{\nu}^t:M \to \mathbb{R}$ defined by 
	\begin{align*}
	 \dist_{\nu}^t(x):=\int _{M} \dist_{M} (x,y)^t \ d \nu (y),
	 \end{align*}where $\dist_M$ is the Riemannian distance on $M$. The proof of the following lemma is easy and we omit the proof.
	\begin{lem}\label{aa}There exists a point $x_{\nu}^t \in M$ such that $\dist_{\nu}$ attains
	 its minimum and  
	 \begin{align*}
	  \dist_{M}(x_{\nu}^t,\supp \nu)\leq 2\diam (\supp \nu).
	  \end{align*}
	 \end{lem}

	\begin{rem}\upshape
	 The above $x_{\nu}^t$ is not unique in general for $0< t
	 \leq  1$. For example,
	 consider $M=\mathbb{R}$ and a Borel probability measure
	 $\nu$ on $\mathbb{R}$ given by $\nu\big(\{0\}\big) = \nu\big(\{
	 1  \}\big)=1/2$. In this case, the function $\dist_{\nu}^t$
	 attains its infimum at both $0$ and $1$. 
	 \end{rem}
	 For $0<s\leq t$, we indicate by $I_M(\nu;s,t)$ the set of
	 all $x\in M$ such that $\dist_{\nu}^r(x)= \min\limits_{y\in
	 M}\dist_{\nu}^r (y)$ for some $s\leq r \leq t$.

	From now on, we consider the $n$-dimensional hyperbolic space $\mathbb{H}^n$ as a Poincar\'{e} disk
	model $D^n:=\{ x\in \mathbb{R}^n  \mid |x|< 1   \}$. For
	$\kappa_1<0$, we denote by $\mathbb{H}^n(\kappa_1)$
	a complete simply connected
	Riemannian manifold of constant sectional curvature $\kappa_1$. We
	consider $(\mathbb{H}^n(\kappa_1),
	\dist_{\mathbb{H}^n(\kappa_1)})$ as $(D^n,
	(1/\sqrt{-\kappa_1})\dist_{\mathbb{H}^n})$.  
	\begin{lem}[{cf.~\cite[Theorem 4.6.1]{rat}}]\label{onemoretime}For any $x,y \in \mathbb{H}^n$ we have
	 \begin{align*}
	  \dist_{\mathbb{H}^n} (x,y)= 2 \log
	  \Big(\frac{|x-y|}{\sqrt{1-|x|^2} \sqrt{1-|y|^2}} + \sqrt{\frac{|x-y|^2}{(1-|x|^2)(1-|y|^2)}+1} \Big). 
	  \end{align*}
	 \end{lem}

For each $n\in \mathbb{N}$, we define the function $\phi_n:D^n \to \mathbb{R}^n$ by 
\begin{align*}
\phi_n (x):=\frac{x}{1-|x|}.
\end{align*}We consider the distance function
$(\phi_n)^{\ast}\dist_{\mathbb{R}^n}$ on $D^n$ defined by
$(\phi_n)^{\ast}\dist_{\mathbb{R}^n}(x,y):= |\phi_n(x)-\phi_n(y)|$. 
\begin{lem}\label{yukiesine}For any $x,y \in D^n$, we have
 $(\phi_n)^{\ast}\dist_{\mathbb{R}^n}(x,y)\geq |x-y|$.
\begin{proof}Observe that $|\widetilde{x}-r\widetilde{y}|\geq |\widetilde{x}-\widetilde{y}|$ for any $\widetilde{x}, \widetilde{y}\in
 \mathbb{R}^n$ and $r\geq 1$. Hence, assuming $|x|\leq |y|$, we obtain
\begin{align*}
(\phi_n)^{\ast}\dist_{\mathbb{R}^n}(x,y)=\Big|\frac{x}{1-|x|}-
 \frac{y}{1-|y|}\Big|=\frac{1}{1-|x|} \Big| x-\frac{1-|x|}{1-|y|}y \Big|
 \geq \frac{|x-y|}{1-|x|} \geq |x-y|.
\end{align*}This completes the proof.
\end{proof}
\end{lem}

Given $\kappa_1<0$, we indicate by $\mathcal{NM}^n(\kappa_1)$
the set of all $n$-dimensional
complete simply connected Riemannian manifolds with nonpositive sectional
curvature bounded below by $\kappa_1$. For $s\geq 0$ and $\kappa_1<0$,
we define $\mathcal{NM}^n(s;\kappa_1)$ by $\mathcal{NM}^n$ if $s=0$
and $\mathcal{NM}^n(\kappa_1)$ if $s>0$.
One of our main theorems in this paper is the following:

\begin{thm}\label{th}
			 Let $\{ (X_n,\dist_n ,
			 \mu_n) \}_{n=1}^{\infty}$ be a sequence of
			 mm-spaces with finite diameter. We assume that
			 sequences $\{a(n) \}_{n=1}^{\infty}$, $\{
 p_n\}_{n=1}^{\infty}$, $\{ s_n\}_{n=1}^{\infty}$,
 $\{\kappa_{n}\}_{n=1}^{\infty}$ of real numbers satisfy that
\begin{align*}
a(n)\in \mathbb{N}, \ p_n >0, \ s_n\geq 0, \ \kappa_n<0, \ 
 \sup\limits_{n\in \mathbb{N}}p_n <+\infty, \ \sup\limits_{n\in
 \mathbb{N}}s_n<1/2 ,\ \sup\limits_{n\in \mathbb{N}} \kappa_n <0, 
\end{align*}and 
			 \begin{align*}
			  \diam \Big(X_{n} \oblip
			  \mathbb{R}^{a(n)}, m_{n}- \frac{\kappa}{
			  (-\kappa_{n})^{s_n}(\diam X_{n})^{p_n} }\Big)
		 \to 0 \text{ as }n\to \infty
			  \end{align*}for any $\kappa>0$. Then for any
			 $\kappa>0$, we have 
	  \begin{align*}
	   \sup \Big\{ \diam\Big(X_{n} \oblip
	   N,m_{n}-\frac{\kappa}{(-\kappa_{n})^{s_n}(\diam X_n )^{p_n}}\Big) \Big| N\in
	   \mathcal{NM}^{a(n)}(s_n;\kappa_n)\Big\}\to 0 \text{ as }n\to \infty.
	   \end{align*}If in addition $\inf\limits_{n\in \mathbb{N}}m_n >0$,
			 then for any sequences $\{ N_n
			 \}_{n=1}^{\infty}$ with $N_n \in
			 \mathcal{NM}^{a(n)}(s_n;\kappa_n)$ and $\{ f_n:X_n
			 \to N_n\}_{n=1}^{\infty}$ of $1$-Lipschitz
			 maps, we have
			 \begin{align*}
			  \diam I_{N_n}\big( (f_n)_{\ast}(\mu_n);\max \{
			  1-2s_n,1-p_n\},1    \big) \to 0 \text{ as }n\to \infty.
			  \end{align*}

		\begin{proof}The claim
		 obviously holds in the case
		 of 
		 $\lim\limits_{n\to \infty}\diam X_n =0$, so we assume that 
		 $\inf\limits_{n\in
		 \mathbb{N}} \diam X_n >0$. 
		 Since
		 $\lim\limits_{n\to
		 \infty}m_{n} =0$ implies the first part of Theorem \ref{th}, we also assume that
		 $\inf\limits_{n \in \mathbb{N}} m_{n} >0$. Note that
		 there is a sequence $\{
		 \widetilde{\kappa}_n\}_{n=1}^{\infty}$ of positive
		 numbers such that $\widetilde{\kappa}_n\to 0$ as $n\to
		 \infty$ and for any $\kappa>0$
		 \begin{align*}
		  \diam \Big( X_n \oblip
		  \mathbb{R}^{a(n)},m_n-\frac{\kappa
		  \widetilde{\kappa}_n}{(-\kappa_n)^{s_n} (\diam
		  X_n)^{p_n}}\Big)\to 0 \text{ as }n\to \infty.
		  \end{align*}Hence we only consider the case that
		 $\kappa_n \to -\infty$ as $n\to \infty$. 
Let $\{  N_n   \}_{n=1}^{\infty}$ be any sequence such that $N_n\in
		 \mathcal{NM}^{a(n)}(s_n;\kappa_n)$ for each $n\in \mathbb{N}$ and $\{ f_n:X_{n}  \to N_n    \}_{n=1}^{\infty}$ 
		 any sequence of $1$-Lipschitz maps. Given
		 arbitrary $t_n$ with $\max \{ 1-2s_n,1-p_n\} \leq t_n \leq
		 1$ and $z_n\in N_n$ with
		 $\dist^{t_n}_{(f_{n})_{\ast}(\mu_n)}(z_n)
		 =\min\limits_{z\in N_n}
		 \dist^{t_n}_{(f_n)_{\ast}(\mu_n)}(z)$, we shall show that
		 \begin{align*}
		  (f_n)_{\ast}(\mu_n) \big( N_n \setminus B_{N_n}(z_n,r)
		  \big) \leq \frac{\kappa }{(-\kappa_{n})^{s_n} (\diam X_n)^{p_n} }
		  \end{align*}for
		 any $r,\kappa>0$ and sufficiently large $n\in \mathbb{N}$. 
		 Take $r_n \geq 1$ with $\overline{f_n(X_n)}\subseteq
		 B_{N_n}(z_n,r_n)$ 
		 and let $\kappa_{1n}$ be a negative number 
		 such that $\kappa_{1n} \to -\infty$ as $n\to \infty$
		 and the
		 sectional curvature on $B_{N_n}(z_n,r_n)$ is bounded
		 from below by $\kappa_{1n}$. Define the function
		 $\kappa_{1n}=\kappa_{1n}(s): [0,+\infty) \to
		 \mathbb{R}$ by $\kappa_{1n}(0):=\kappa_{1n}$ for $s=0$ and $\kappa_{1n}(s):=
		 \kappa_n$ for $s>0$. We observe that $\kappa_{1n}(s_n)\to
		 -\infty$ as $n\to \infty$. 
		 Let
		 $\widetilde{\varphi}_n$ be a linear isometry from the
		 tangent space of $N_n$ at $z_n$ to the tangent space of $
		 \mathbb{R}^{a(n)}$ at $0$ and put
		 $\varphi_n:=\phi_{a(n)}^{-1} \circ \exp_{0} \circ \widetilde{\varphi}_n \circ
		 \exp_{z_n}^{-1}:N_n \to D^{a(n)}$. By virtue of the hinge theorem
 (cf.~{\cite[Chapter \Roman{yon}, Remark $2.6$]{sakai}}) and Lemma \ref{yukiesine}, we have
\begin{align*}
  |\varphi_n(x)-\varphi_n(x')|   \leq (\phi_{a(n)})^{\ast}\dist_{\mathbb{R}^{a(n)}}\big(
 \varphi_n(x),\varphi_n(x') \big)\leq \dist_{N_{n}}(x,x') 
\end{align*}for any $x,x'\in N_n$. Since ${\varphi}_n \circ f_n$
		  is the $1$-Lipschitz map from $X_{n}$ to the Euclidean space $(D^{a(n)},\dist_{\mathbb{R}^{a(n)}})$, from the assumption there
		 exists a Borel
		 subset $A_n \subseteq D^{a(n)}$ and $\widetilde{\kappa}_n >0$
		 such that
		 $A_n \subseteq \supp ({\varphi}_n \circ f_n
		 )_{\ast}(\mu_{n})\subseteq \varphi_n\big( B_{N_n}(z_n ,
		 r_n)\big)$, $({\varphi}_n \circ
		  f_n)_{\ast}(\mu_{n})(A_n)> m_{n}-\widetilde{\kappa}_n
		 /\big((-\kappa_{n})^{s_n}(\diam X_n)^{p_n}\big)$ and $\diam
		 (A_n,\dist_{\mathbb{R}^{a(n)}}) \to 0 $ as $n\to \infty$.

		 Let us show that $
		  (\phi_{a(n)})^{\ast}\dist_{\mathbb{R}^{a(n)}}
		 (A_n,0)\to 0$ as $n\to
		 \infty$. Suppose that there exists a constant $C>0$
		 such that $(\phi_{a(n)})^{\ast}\dist_{\mathbb{R}^{a(n)}} (A_n
		  ,0) \geq C$ for any $n\in \mathbb{N}$. For any $x\in D^{a(n)}$, we take
 $y\in D^{a(n)}$ such that $y=\lambda x$ as a vector in
		 $\mathbb{R}^{a(n)}$ for some $\lambda \geq 0$ and 
\begin{align*}
\frac{1}{\sqrt{-\kappa_{1n}(s_n)}}\dist_{\mathbb{H}^{a(n)}}(0,y)=(\phi_{a(n)})^{\ast}\dist_{\mathbb{R}^{a(n)}}(0,x).
\end{align*}We define $\psi_n :D^{a(n)}\to D^{a(n)}$ by $\psi_n
 (x):=y$.  
		 Since $\diam (A_n, \dist_{\mathbb{R}^{a(n)}}) \to 0$, $\kappa_{1n}(s_n) \to -\infty$ as $n\to
		 \infty$, and $(\phi_{a(n)})^{\ast}\dist_{\mathbb{R}^{a(n)}}(A_n,0)\geq
		 C$, we get 
		 \begin{align*}
		  \diam(\psi_n(A_n),\dist_{\mathbb{R}^{a(n)}}), \dist_{\mathbb{R}^{a(n)}}(\psi_n(A_n),\mathbb{S}^{a(n)-1})
		  \to 0 \text{ as }n\to \infty. 
		  \end{align*}Since $\diam
		 (\psi_n(A_n),\dist_{\mathbb{R}^{a(n)}})$, $\dist_{\mathbb{R}^{a(n)}}(\psi_n(A_n) ,
		 \mathbb{S}^{a(n)-1})\to 0$ as $n\to \infty$, there are
		 points $q_n\in
		 (\psi_n \circ \varphi_n )\big(B_{N_n}(z_n,r_n)\big)=
		 B_{\mathbb{H}^{a(n)}}(0,\sqrt{-\kappa_{1n}(s_n)}r_n)$ having the following
		 properties $(1)$, $(2)$: 
\begin{itemize}
\item[$(1)$]$|q_n|\to 1$ as $n\to \infty$.
\item[$(2)$]$|q_n-x|\leq 2(1-|q_n|)$ and $|q_n|\leq |x|$ for any $x\in
	    \psi_n (A_n)$.
\end{itemize}

\begin{claim}\label{karee}For any $n\in \mathbb{N}$ and $x\in
		  \varphi_n^{-1}(A_n)\cap f_n(X_n)$, we have 
\begin{align*}
\dist_{N_n}(z_n,x)^{t_n}\geq
\dist_{N_n}\big((\psi_n \circ \varphi_n)^{-1}(q_n),x\big)^{t_n}+\frac{c_n}{(\diam X_n)^{1-t_n}},
\end{align*}where
 $c_n$ are some positive numbers satisfying $\inf\limits_{l\in \mathbb{N}}\big(c_l/\dist_{N_l}\big(z_l,(\psi_l \circ \varphi_l)^{-1}(q_l)\big)\big) >0$. 
\begin{proof}
According to the hinge theorem (cf.~{\cite[Chapter \Roman{yon}, Theorem
 $4.2$ $(2)$]{sakai}}), we have $\dist_{N_n}\big((\psi_n \circ
 \varphi_n)^{-1}(q_n),x \big)\leq
 (1/\sqrt{-\kappa_{1n}(s_n)})\dist_{\mathbb{H}^{a(n)}}\big(  q_n,
 (\psi_n\circ \varphi_n)(x)       \big)$. 
Note that $\dist_{N_n}(z_n,x)=(1/\sqrt{-\kappa_{1n}(s_n)})\dist_{\mathbb{H}^{a(n)}}\big(0,
 (\psi_n \circ \varphi_n)(x)\big)$. Therefore, 
from Lemma \ref{onemoretime} and $(2)$, we have

		  \begin{align*}
&\sqrt{-\kappa_{1n}(s_n)}\dist_{N_n}(z_n,x)-
		   \sqrt{-\kappa_{1n}(s_n)}\dist_{N_n}\big( (\psi_n \circ \varphi_n)^{-1}(q_n),x\big)\\
\geq \  
&  \dist_{\mathbb{H}^{a(n)}}\big(0,(\psi_n
		   \circ \varphi_n)(x)\big)-
		   \dist_{\mathbb{H}^{a(n)}}\big(q_n , (\psi_n \circ
		   \varphi_n )(x)\big)\\
= \  & 2\log \Big(\frac{\sqrt{1-|q_n|^2}(1+|(\psi_n \circ
		   \varphi_n)(x)| )    }{|q_n -(\psi_n \circ
		   \varphi_n)(x)|+ \sqrt{|q_n - (\psi_n \circ
		   \varphi_n)(x)|^2  +(1-|(\psi_n \circ
		   \varphi_n)(x)|^2)(1-|q_n|^2)          }   }\Big) \\
\geq \ & 2\log \frac{1}{2}\Big(\frac{\sqrt{1-|q_n|^2}(1+|(\psi_n \circ
		   \varphi_n)(x)| )    }{1-|q_n|+ \sqrt{(1-|q_n|)^2  +(1-|(\psi_n \circ
		   \varphi_n)(x)|^2)(1-|q_n|^2)          }   }\Big) \\
= \ & 2\log \frac{1}{2} \Big(      \frac{(1+|q_n|)(1+|(\psi_n \circ \varphi_n)(x)|) }{\sqrt{1-|q_n|^2}+ \sqrt{1-|q_n|^2  +(1-|(\psi_n \circ
		   \varphi_n)(x)|^2)(1+|q_n|)^2          }   }
		   \Big) \\ 
\geq \ &2 \log \frac{1}{2}  \Big(      \frac{(1+|q_n|)^2
		   }{\sqrt{1-|q_n|^2} \big(1+  \sqrt{1+ (1+|q_n|)^2}    \big)  }
		   \Big)=:\sqrt{-\kappa_{1n}(s_n)}b_n.
		   \end{align*}Thus, combining this with $t_n\leq 1$, we get
\begin{align*}
& \dist_{N_n}(z_n,x)^{t_n}- \dist_{N_n}\big((\psi_n \circ
 \varphi_n)^{-1}(q_n),x \big)^{t_n}\\ 
\geq
 \ &
 \frac{\dist_{N_n}(z_n,x)^{t_n}-\dist_{N_n}\big( (\psi_n \circ
 \varphi_n)^{-1}(q_n), x\big)^{t_n}}{\dist_{N_n}(z_n,x)-\dist_{N_n}\big( (\psi_n \circ
 \varphi_n)^{-1}(q_n), x\big)}\big(  \dist_{N_n}(z_n,x)-\dist_{N_n}\big( (\psi_n \circ
 \varphi_n)^{-1}(q_n), x\big)\big)
 \\ \geq \ &\frac{b_n t_n}{\big(\theta \dist_{N_n}(z_n,x)+(1-\theta)\dist_{N_n}\big( (\psi_n \circ
 \varphi_n)^{-1}(q_n), x\big)\big)^{1-t_n}} \ \ (\text{ for some } 0\leq
 \theta \leq 1)\\
 \geq \ &
 \frac{b_n t_n}{\dist_{N_n}(z_n,x)^{1-t_n}}.
\end{align*}Applying Lemma \ref{aa}, we have $\dist_{N_n}(z_n,x)\leq 2\diam
 X_n$. Therefore we obtain
\begin{align*}
\dist_{N_n}(z_n,x)^{t_n} - \dist_{N_n}\big((\psi_n\circ \varphi_n)^{-1}(q_n),x)^{t_n}
 \geq \frac{b_n t_n}{2^{1-t_n} (\diam X_n)^{1-t_n}}.
\end{align*}Since $ b_n /
 \dist_{N_n}\big(z_n,(\psi_n \circ \varphi_n)^{-1}(q_n)\big) \to 1$ as
 $n\to \infty$ and $0<\inf\limits_{l\in \mathbb{N}}(1-2s_l) \leq t_n
 \leq 1$, putting $c_n:=2^{t_n-1}b_n t_n $, this completes the proof of the claim.
\end{proof}
			 \end{claim}

Put 
\begin{align*}
B_n : = \{ x_n \in X_n \setminus (\varphi_n \circ
 f_n)^{-1}(A_n) \mid \dist_{N_n}\big( (\psi_n \circ
		 \varphi_n)^{-1}(q_n), f_n(x_n)   \big) \geq
 \dist_{N_n}\big(z_n,f_n(x_n)\big)\}. 
\end{align*}By virtue of Claim \ref{karee}, we have
		  \begin{align}\label{nekounko2}
	   &\dist^{t_n}_{(f_{n})_{\ast}(\mu_n)}\big((\psi_n  \circ \varphi_n)^{-1}(q_n)\big)
\\ 
		   = &\int_{(\varphi_n \circ f_n)^{-1}(A_n)}
		   \dist_{N_n}\big((\psi_n \circ \varphi_n)^{-1}(q_n),  f_n(x_n)    \big)^{t_n}
		   \ d \mu_{n} (x_n ) \tag*{} \\
		   &\hspace{5cm} + \int_{X_n \setminus (\varphi_n \circ f_n)^{-1}(A_n)} \dist_{N_n}\big((\psi_n
		   \circ \varphi_n)^{-1}(q_n) , f_n(x_{n})\big)^{t_n} \ d \mu_{n}
		   (x_{n}) \tag*{} \\
		   \leq  & \int_{({\varphi}_n \circ f_n)^{-1}(A_n)} \dist_{N_n}\big(z_n,f_n(x_{n})\big)^{t_n} \ d \mu_{n} (x_{n})
		   - \frac{c_n}{(\diam X_n)^{1-t_n}} ({\varphi}_n \circ
		   f_n)_{\ast}(\mu_{n})(A_n) \tag*{} \\
		    & \hspace{5cm} + \int_{X_n\setminus (\varphi_n \circ
		   f_n)^{-1}(A_n)} \dist_{N_n}\big((\psi_n
		   \circ \varphi_n)^{-1}(q_n) , f_n(x_{n}) \big)^{t_n} \ d
		   \mu_{n} (x_n) \tag*{} \\
		   \leq & \dist_{(f_{n})_{\ast}(\mu_n)}^{t_n}(z_n)
		   -\frac{c_n}{(\diam X_n)^{1-t_n}}\Big(\inf_{l \in
		   \mathbb{N}}m_{l}-\frac{\widetilde{\kappa}_n}{(-\kappa_{n})^{s_n}(\diam
		   X_n)^{p_n}}\Big) \tag*{} \\
		    & \hspace{3cm} + \int_{B_n} \big\{ \dist_{N_n}\big(
		   (\psi_n \circ \varphi_n)^{-1}(q_n),f_n(x_n)
		   \big)^{t_n}-\dist_{N_n}\big(z_n,f_n(x_n)\big)^{t_n}
		   \big\} \ d
		   \mu_{n} (x_{n}) \tag*{} \\
\leq & \dist_{(f_{n})_{\ast}(\mu_n)}^{t_n}(z_n)
		   -\frac{c_n \inf\limits_{l \in
		   \mathbb{N}}m_{l}}{2(\diam X_n)^{1-t_n}} \tag*{} \\
		    & \hspace{4cm} 
+ \int_{B_n} \frac{ \dist_{N_n}\big(
		   (\psi_n \circ \varphi_n)^{-1}(q_n),f_n(x_n)
		   \big)-\dist_{N_n}\big(z_n,f_n(x_n) \big) }{ \dist_{N_n}\big(
		   (\psi_n \circ \varphi_n)^{-1}(q_n),f_n(x_n)
		   \big)^{1-t_n}}
		    \ d
		   \mu_{n} (x_{n}) \tag*{}
		   \end{align}for any sufficiently large $n\in
		 \mathbb{N}$. 
	
		 Assume first that $\mu_n(B_n)=0$ for infinitely many $n\in
		 \mathbb{N}$. From the above inequality, we get
		 \begin{align*}\dist^{t_n}_{(f_n)_{\ast}(\mu_n)} \big(
		  (\psi_n \circ \varphi_n)^{-1}(q_n) \big) \leq
		  \dist^{t_n}_{(f_n)_{\ast}(\mu_n)} (z_n) - \frac{c_n
		  \inf\limits_{l\in \mathbb{N}}m_l}{2(\diam
		  X_n)^{1-t_n}}< \dist_{(f_n)_{\ast}(\mu_n)}^{t_n}(z_n).
		  \end{align*}for any sufficiently large $n\in
		 \mathbb{N}$. This is a contradiction since $z_n \in
		 N_n$ is the infimum of the function
		 $\dist^{t_n}_{(f_n)_{\ast}(\mu_n)}$.

		 We consider the other case that $\mu_n(B_n)>0 $ for any sufficiently large
		 $n\in \mathbb{N}$. Since $2\dist_{N_n}\big((\psi_n \circ
		 \varphi_n)^{-1}(q_n),f_n(x_n)\big) \geq \dist_{N_n}\big(z_n,(\psi_n \circ
		 \varphi_n)^{-1}(q_n)\big)$ for any $x_n \in B_n$ and
		 $s_n\geq 1-t_n$, we have 
		 \begin{align}\label{nekounko}
		  &\Big(\frac{-\kappa_{n}}{\inf\limits_{l \in
		  \mathbb{N}}(-\kappa_{l})}\Big)^{s_n}
		  \dist_{N_n}\big((\psi_n \circ \varphi_n)^{-1}(q_n),f_n(x_n)\big)^{1-t_n}  \\
		  \geq & \ \big(\inf\limits_{l\in \mathbb{N}}(-\kappa_l)\big)^{\frac{t_n-1}{2}}
		\big(\sqrt{-\kappa_{1n}(s_n)}\dist_{N_n}\big(z_n,  (\psi_n \circ \varphi_n)^{-1}(q_n)
		 \big)/2 \big)^{1-t_n} \tag*{} \\ 
		  = & \ \big(\inf\limits_{l\in \mathbb{N}}(-\kappa_l)\big)^{\frac{t_n-1}{2}}\big(\dist_{\mathbb{H}^{a(n)}}( 0,   q_n
		  )/2 \big)^{1-t_n} \geq 1 \tag*{}
		  \end{align}for any sufficiently large $n\in
		 \mathbb{N}$ and $x_n \in B_n$. Since $p_n\geq
		 1-t_n$, $\sup\limits_{l\in \mathbb{N}}p_l <+\infty$, 
		 and $\inf\limits_{l\in \mathbb{N}}\diam X_l >0$, we
		 obtain $\inf\limits_{l\in \mathbb{N}}(\diam
		  X_l)^{p_l + t_l -1}>0$. Combining this with (\ref{nekounko}), $\widetilde{\kappa}_n \to 0$ as $n\to \infty$, 
		 $\inf\limits_{l\in
		 \mathbb{N}}\big(c_l/\dist_{N_l}\big(z_l,(\psi_l \circ
		 \varphi_l)^{-1}(q_l)\big)\big) >0$, and $\mu_n(B_n)<
		 \widetilde{\kappa}_n/ \big( (-\kappa_n)^{s_n}(\diam X_n)^{p_n}\big)$, we get
		  \begin{align}\label{nekounko3}
 &-\frac{c_n \inf\limits_{l \in
		   \mathbb{N}}m_{l}}{2(\diam X_n)^{1-t_n}} 
+ \int_{B_n} \frac{ \dist_{N_n}\big(
		   (\psi_n \circ \varphi_n)^{-1}(q_n),f_n(x_n)
		   \big)-\dist_{N_n}\big(z_n,f_n(x_n) \big) }{ \dist_{N_n}\big(
		   (\psi_n \circ \varphi_n)^{-1}(q_n),f_n(x_n)
		   \big)^{1-t_n}}
		    \ d
		   \mu_{n} (x_{n})\\
\leq \ & \int_{B_n} \Big\{\frac{\dist_{N_n}\big(z_n,
		   (\psi_n \circ \varphi_n)^{-1}(q_n)
		   \big) }{\dist_{N_n}\big(
		   (\psi_n \circ \varphi_n)^{-1}(q_n),f_n(x_n)
		   \big)^{1-t_n}   }- \frac{c_n (-\kappa_n)^{s_n}(\diam X_n)^{p_n+t_n-1}\inf\limits_{l\in
		   \mathbb{N}}m_l}{2\widetilde{\kappa}_n}\Big\} \ d
		   \mu_{n} (x_{n}) \tag*{} \\
  \leq \ & \int_{B_n} \frac{\dist_{N_n}\big(z_n,
		   (\psi_n \circ \varphi_n)^{-1}(q_n)
		   \big)-2\dist_{N_n}\big( z_n, (\psi_n \circ \varphi_n)^{-1}(q_n)\big)}{\dist_{N_n}\big(
		   (\psi_n \circ \varphi_n)^{-1}(q_n),f_n(x_n) \big)^{1-t_n}       } \ d
		   \mu_{n} (x_{n}) <0 \tag*{}
		   \end{align}for any sufficiently large $n\in
		 \mathbb{N}$. Therefore, from (\ref{nekounko2}) and
		 (\ref{nekounko3}), we obtain 
\begin{align*}
\dist_{(f_{n})_{\ast}(\mu_{n})}^{t_n}\big(
		 (\psi_n \circ \varphi_n)^{-1}(q_n)\big) < \dist_{(f_{n}
		 )_{\ast}(\mu_{n})}^{t_n}(z_n)
\end{align*}
for any sufficiently large $n\in
		 \mathbb{N}$. Consequently we have a
		  contradiction since $z_n\in N_n$ is the infimum of the function $\dist_{(f_{n})_{\ast}(\mu_{n})}^{t_n}$.

	Since $(\phi_{a(n)})^{\ast}\dist_{\mathbb{R}^{a(n)}}(A_n,0)\to
		 0$ as $n\to \infty$, we get $A_n \subseteq \varphi_n \big(B_{N_n}(z_n,
		 r)\big)$ for any $r>0$ and sufficiently large $n\in
		 \mathbb{N}$. Therefore we obtain 
\begin{align*}
		  (f_n)_{\ast}(\mu_n) \big( N_n \setminus B_{N_n}(z_n,r)
		  \big) \leq (f_n)_{\ast}(\mu_n) \big( N_n \setminus \varphi_n^{-1}(A_n)
		  \big) \leq \frac{\kappa }{(-\kappa_{n})^{s_n}(\diam X_n)^{p_n}}
\end{align*}for any $\kappa>0$ and sufficiently large $n\in \mathbb{N}$. This
		 completes the proof of the theorem.
	\end{proof}
\end{thm}

\begin{proof}[Proof of Theorem \ref{kuzuninaritaiwan}]Let $\{
 N_n\}_{n=1}^{\infty}$ be any sequence such that $N_n\in
 \mathcal{NM}^{a(n)}$ and $\{f_n :X_n\to
 N_n\}_{n=1}^{\infty}$ any sequence of
 $1$-Lipschitz maps. We shall show that 
 \begin{align}\label{seikeibijo}
\diam\big((f_n)_{\ast}(\mu_n),m_n-\kappa\big)\to 0\text{ as }n\to \infty
\end{align}for any $\kappa>0$. By virtue of Corollary \ref{kawaisouni2},
 there is a Borel subset $A_n
 \subseteq X_n$ such that $\diam A_n <+\infty$, $m_n':=\mu_n(A_n)\geq m_n
 -\kappa /2$, and 
\begin{align}\label{kawaisouni3}
\diam (A_n \oblip \mathbb{R}^{a(n)},m_n'-\kappa')\to 0 \text{ as }n\to \infty
\end{align}for any $\kappa'>0$. The claim (\ref{seikeibijo}) obviously holds in the case
 of $\lim\limits_{n\to \infty}\diam A_n =0$, so we assume that
 $\inf\limits_{n\in \mathbb{N}}\diam A_n>0$. Observe that 
\begin{align*}
0< \inf_{n\in \mathbb{N}}(\diam A_n)^{1/\diam A_n} \leq \sup_{n\in
 \mathbb{N}}(\diam A_n)^{1/\diam A_n}<+\infty.
\end{align*}From this and (\ref{kawaisouni3}), we have 
\begin{align*}
\diam \Big( A_n \oblip \mathbb{R}^{a(n)}, m_n'-\frac{\kappa'}{(\diam
 A_n)^{1/\diam A_n}}\Big) \to 0 \text{ as }n\to \infty
\end{align*}for any $\kappa'>0$. Combining this with the same proof in Theorem \ref{th} for
 the sequence $\{f_n|_{A_n} :A_n \to N_n\}_{n=1}^{\infty}$, we obtain 
\begin{align*}
\diam((f_n)_{\ast}(\mu_n),m_n-\kappa)\leq \diam \big((f_n|_{A_n})_{\ast}(\mu_n),m_n'-\kappa/2 \big) \to 0 \text{ as }n\to \infty.
\end{align*}This completes the proof.
\end{proof}

\begin{proof}[Proof of Theorem \ref{koukaithm}]We first assume that $f_{\ast}(\mu_X)
 \in B^1(N)$ for any $1$-Lipschitz map $f:X\to N$. Given an arbitrary
 $1$-Lipschitz map $f:X\to N$, we put $z:=b\big(f_{\ast}(\mu_X)\big)$. From the triangle inequality, we have
\begin{align}\label{423}
 V_p(f)\leq  2 \Big(\int\int_{N\times N} \dist_N(x,z)^p\
 df_{\ast}(\mu_X)(x) df_{\ast}(\mu_X)(y)\Big)^{1/p}
 =  2 \Big(m\int_{N} \dist_N(x,z)^p\
 df_{\ast}(\mu_X)(x)\Big)^{1/p}.
 \end{align}We identify the tangent space of $N$ at $z$ with the
 Euclidean space $\mathbb{R}^n$ and consider the map $f_0:=\exp^{-1}_z
 \circ f:X\to \mathbb{R}^n$. According to the hinge theorem, the map $f_0$ is
  $1$-Lipschitz. Since the map $\exp^{-1}_z$ is isometric on rays
 issuing from $z$, we get
 \begin{align}\label{424}
  \int_{N} \dist_N(x,z)^p\
 df_{\ast}(\mu_X)(x) = \int_N |\exp^{-1}_zx|^p \ df_{\ast}(\mu_X)(x) =
  \int_{\mathbb{R}^n} |y|^p \ d(f_{0})_{\ast}(\mu_X)(y).
  \end{align}Since $b \big((f_0)_{\ast}(\mu_X)\big) =0$ by Lemma \ref{atatata}, it
 follows from Lemma \ref{atata} that 
 \begin{align}\label{425}
  |y|^p = \Big| \frac{1}{m}\int_X (y-y')\ d(f_0)_{\ast}(\mu_X)(y')\Big|^p\leq \frac{1}{m}\int_{\mathbb{R}^n} |y-y'|^p\ d(f_{0})_{\ast}(\mu_X)(y').
  \end{align}Therefore, combining (\ref{423}) with (\ref{424}) and (\ref{425}), we obtain $V_p
 (f) \leq 2V_p(f_0)$.
 
 We consider the other case that there exist a $1$-Lipschitz map $f:X\to
 N$ with $f_{\ast}(\mu_X)\notin B^1(N)$. From H\"{o}lder's inequality,
 we have $\int_N \dist_N (x,y)^p\ df_{\ast}(\mu_X)(y)=+\infty$ for any
 $x\in N$. Hence, Fubini's theorem yields $V_p(f)=+\infty$. Take a
 point $x_0 \in X$. For each $k\in \mathbb{N}$, we put
 $f_k:=f|_{B_X(x_0,k)}$. Since $(f_k)_{\ast}(\mu_X)\in B^1(N)$, from the
 above proof, there exists a $1$-Lipschitz map $\widetilde{f}_k :
 B_X(x_0,k)\to \mathbb{R}^n$ such that $V_p (f_k)\leq 2V_p
 (\widetilde{f}_k)$. From \cite[Theorem 3.1.2]{ambro}, there exists $\sqrt{n}$-Lipschitz extension of
 $\widetilde{f}_k$, say $g_k:X\to \mathbb{R}^n$. Since
 $(1/\sqrt{n})g_k:X\to \mathbb{R}^n$ is a $1$-Lipschitz map and
\begin{align*}
V_p\Big(\frac{1}{\sqrt{n}}g_k\Big)  \geq
 \frac{1}{\sqrt{n}}V_p(\widetilde{f_k}) \geq
 \frac{1}{2\sqrt{n}}V_p(f_k)\to \frac{1}{2\sqrt{n}}V_p(f)=+\infty \text{ as }k\to \infty, 
\end{align*}we obtain $\obin_{\mathbb{R}^n}(X)=+\infty$. This completes
 the proof of the theorem.
 \end{proof}

 \begin{rem}\label{katoukatou}\upshape Let us consider the case of $p=2$. Since
  \begin{align*}
   \int_{\mathbb{R}^n} |y|^2 \ d(f_0)_{\ast}(\mu_X)(y) =
   \frac{1}{2m}\int\int_{\mathbb{R}^n \times \mathbb{R}^n} |y-y'|^2 \
   d(f_0)_{\ast}(\mu_X) (y) d(f_0)_{\ast}(\mu_X)(y')
   \end{align*}for any $1$-Lipschitz maps $f_0:X\to \mathbb{R}^n$ with
  mean zero, a slight modification of the proof of Theorem
  \ref{koukaithm} implies $\obinin_N(X)\leq \sqrt{2}\obinin_{\mathbb{R}^n}(X)$. 
  \end{rem}

\begin{proof}[Proof of Theorem \ref{kameda1}]Applying Theorem \ref{th}
 to Corollary \ref{waruiyokan}, we obtain the proof of the theorem. 
\end{proof}
	 \begin{proof}[Proof of Proposition \ref{kantan}]By $(2)$ and $m_n < +\infty$, $X_n$ is 
	  compact. Hence, there exist $x_n, y_n \in X_n$ such that $\diam X_n = \dist_n (x_n,y_n)$. We define a function $f_n:X_n \to
	  \mathbb{R}$ by $f_n(x_n):=\dist_n(x_n,x)$. Let $a_n$ be a pre-L\'{e}vy mean of $f_n$. Then we have 
	  \begin{claim}\label{cl1}$ \diam X_n /2 \leq a_n$.
	   \begin{proof}If $a_n < \diam X_n/2$, we get $B_{X_n}(x_n,a_n)\cap B_{X_n}(y_n,a_n)=\emptyset$. Since $B_{X_n}(x_n,a_n)=\{
		x\in X_n \mid f_n(x)\leq a_n     \}$ and by $(2)$, we obtain 
	  \begin{align*}
	   m_n \leq 2\mu_n\big( B_{X_n}(x_n,a_n) \big)=\ &\mu_n\big(B_{X_n}(x_n,a_n)\big)+\mu_n\big(  B_{X_n}(y_n,a_n)    \big)\\=\ &\mu_n\big(
	   B_{X_n}(x_n,a_n)\cup B_{X_n}(y_n,a_n) \big),
	   \end{align*}which implies that $X_n=\supp \mu_n \subseteq B_{X_n}(x_n,a_n)\cup B_{X_n}(y_n,a_n)$. This is a contradiction
		because $X_n$ is connected. This completes the proof of the claim.
	   \end{proof}
	   \end{claim}
	  Given $y\in X_n$, we define a function $f_y:X_n \to \mathbb{R}$ by 
	  $f_{y}(x):= \dist_n(x,y)$. Then we have
	  \begin{claim}\label{cl2}$a_n$ is a pre-L\'{e}vy mean of $f_y$.
	   \begin{proof}Since 
	  \begin{align*}
	   \mu_n(\{ x\in X_n \mid f_y(x)\leq a_n    \})=\ &\mu_n\big(B_{X_n}(y,a_n)\big)
	   \\ = \ &\mu_n\big(B_{X_n}(x_n,a_n)\big)=\mu_n(\{ x\in X_n \mid
	   f_n(x)\leq a_n \}),
	   \end{align*}we obtain $\mu_n(\{ x\in X_n \mid f_{y}(x)\leq a_n    \})\geq m_n /2$. On the other hand, we get 
	  \begin{align*}
	   \mu_n(  \{  x\in X_n \mid f_{y}(x)\geq a_n     \})= \ &\lim_{\varepsilon \to +0}\mu_n\big(X_n \setminus B_{X_n}(y,a_n-\varepsilon)\big)\\
	   = \ &\lim_{\varepsilon \to +0}\mu_n\big(X_n \setminus B_{X_n}(x_n,a_n -\varepsilon)\big)\\
	   = \ & \mu_n(\{ x\in X_n \mid f_n(x)\geq a_n      \})\geq \frac{m_n}{2}.
	   \end{align*}This completes the proof.
		\end{proof}
	   \end{claim}
	  Take any $\varepsilon > 0$ and put 
	  $B_{n,\varepsilon}:=\{ (x,y) \in X_n \times X_n \mid |\dist_n (x,y)-a_n|<\varepsilon\}$. Then, it follows from Fubini's
	  theorem together with Claim \ref{cl1} and \ref{cl2} that 
	  \begin{align*}
	  ( \mu_n \times \mu_n) (B_{n,\varepsilon})=\int_{X_n} \mu_n( \{ y\in X_n \mid |f_y(x)-a_n|<\varepsilon     \}) \
	   d\mu_n(y)\geq m_n (m_n -2\alpha_{X_n}(\varepsilon)).
	   \end{align*}Let $\iota_n:X_n\to Y_n$ be an isometric embedding and $A_n \subseteq Y_n$ any Borel subset with $(\iota_{n})
	  _{\ast}(\mu_n)(A_n)\geq m_n -\kappa$. Then we have
	  \begin{align*}
	   &(\mu_n \times \mu_n)\big(\big(\iota_n^{-1}(A_n) \times \iota_n^{-1}(A_n)\big) \cap B_{n,\varepsilon}\big)\\
	   \geq \ &\big((\iota_{n})_{\ast}(\mu_n) \times (\iota_{n} )_{\ast}(\mu_n)\big)(A_n \times A_n)+ (\mu_n
	   \times \mu_n)(B_{n,\varepsilon})-m_n^2
	   \geq (m_n-\kappa)^2 -2m_n \alpha_{X_n}(\varepsilon).
	   \end{align*}Since $\alpha_{X_n}(\varepsilon)\to 0$ as $n\to \infty$, the right-hand side of the
	  above inequality is positive for any sufficiently large $n\in \mathbb{N}$. So, we get $B_{n,\varepsilon}\cap
	  \big(\iota_{n}^{-1}(A_n) \times \iota_{n}^{-1}(A_n)\big)\neq \emptyset$, which leads to that there exist
	  $x_n^0,y_n^0 \in \iota_n^{-1}(A_n)$ with $|\dist_n(x_n^0,y_n^0)-a_n|<\varepsilon$. 
	  Hence, we obtain
	  \begin{align*}
	   \diam (X_n \oblip Y_n,m_n -\kappa)\geq \diam \big((\iota_{n})_{\ast}(\mu_n),m_n -\kappa\big)\geq a_n -\varepsilon \geq \frac{\diam X_n}{2}-\varepsilon,
	   \end{align*}and this completes the proof of the proposition.
	  \end{proof}

 \section{Applications and some related topics about the main theorems}
\subsection{Observable central radius}\label{kondokosou}

Let $Y$ be a metric space and assume that a measure $\nu \in B^1(Y)$ has the
barycenter. For any $\kappa >0$, putting $m:= \nu (Y)$, we define the \emph{central radius}
$\crad(\nu,m-\kappa)$ of $\nu$ as the infimum of $\rho >0$ such that
$\nu\big(B_Y(b(\nu),\rho)\big)\geq m-\kappa$.

Let $(X,\dist_X, \mu_X)$ be an mm-space with $\mu_X \in B^1(X)$ and $Y$
a barycentric metric space. For any $\kappa >0$, we define 
\begin{align*}
\obsc_Y(X;-\kappa):= \sup \{ \crad(f_{\ast}(\mu_X),m-\kappa) \mid f:X\to
 Y \text{ is a }1 \text{-Lipschitz map}\},
\end{align*}and call it the \emph{observable central radius} of $X$.

The proof of the following lemma is easy, so we omit the proof.
\begin{lem}[{cf.~\cite[Section $3\frac{1}{2}.31$]{gromov}}]\label{hayakusitekure}For any $\kappa >0$, we have
\begin{align*}
\diam (X \oblip Y,m-\kappa )\leq 2\obsc_Y(X;-\kappa).
\end{align*}
\end{lem}
We prove the following lemma by Gromov, since we find no proof anywhere
in \cite{gromov}.
\begin{lem}[{cf.~\cite[Section $3\frac{1}{2}.31$]{gromov}}]\label{eminobaka}Let $X$ be an mm-space with a finite diameter and $(Y,\| \cdot \|)$ a Hilbert
 space. Then for any $\kappa >0$ we have
\begin{align}\label{doudemoyoshi}
\obsc_Y(X;-\kappa )\leq \diam(X\oblip Y,m-\kappa) + \frac{\kappa}{m}\diam X.
\end{align}
\begin{proof}Given a $1$-Lipschitz map $f:X\to Y$, let $Y_0 \subseteq
 Y$ be a Borel subset such that $f_{\ast}(\mu_X)(Y_0)\geq m-\kappa$. From Lemma \ref{atata}, we get
 \begin{align*}
	   \| b\big( f_{\ast}(\mu_X) \big)-z \|= \ &\Big\| \frac{1}{m}\int_Y(y-z)\ df_{\ast}(\mu_X)(y)  \Big\|\\
	   \leq \ &\frac{1}{m}\int_Y \| y-z\| \
	   df_{\ast}(\mu_X)(y)\\
	   =\ &\frac{1}{m}\int_{Y_0} \| y-z\| \ df_{\ast}(\mu_X)(y)+\frac{1}{m}\int_{Y\setminus Y_0} \| y-z \| \ df_{\ast}(\mu_X)(y)\\
	   \leq \ &\diam
	   Y_0 +\frac{\kappa}{m}\diam X. 
	   \end{align*}Hence we obtain
 \begin{align*}
	   f_{\ast}(\mu_X)\Big(   B_Y\Big(b\big(f_{\ast}(\mu_X)\big),\diam Y_0 +\frac{\kappa}{m}\diam X \Big)\Big)\geq
	   f_{\ast}(\mu)(Y_0)\geq m-\kappa,
 \end{align*}which implies that $\crad \big(
 f_{\ast}(\mu_X),m-\kappa\big)\leq \diam Y_0 +\frac{\kappa}{m}\diam
 X$. This completes the proof.
\end{proof}
\end{lem}

\begin{rem}\upshape In \cite{gromov}, Lemma \ref{eminobaka} is stated as 
\begin{align}\label{doudemoyoshi2}
\obsc_Y(X;-\kappa) \leq \diam (X\oblip
 Y,m-\kappa)+\frac{\kappa}{m-\kappa} \diam X.
\end{align}The inequality (\ref{doudemoyoshi}) gives slightly better
 estimate than that (\ref{doudemoyoshi2}) gives.
\end{rem}

\begin{lem}\label{morimadoka}For any $N\in \mathcal{NM}^n$ and $\kappa >0$, we have
\begin{align*}
\obsc_N(X;-\kappa)\leq \obsc_{\mathbb{R}^n}(X;-\kappa).
\end{align*}In particular, if sequences $\{X_n\}_{n=1}^{\infty}$ of
 mm-spaces with $\mu_n \in B^1(X_n)$ and $\{ a(n) \}_{n=1}^{\infty}$ of natural numbers
 satisfy that
\begin{align*}
\obsc_{\mathbb{R}^{a(n)}}(X_n;-\kappa)\to 0 \text{ as }n\to \infty 
\end{align*}for any $\kappa >0$, then for any $\kappa>0$ we have
\begin{align*}
\sup \{   \obsc_{N}(X_n;-\kappa) \mid N \in \mathcal{NM}^{a(n)}    \} \to 0
 \text{ as }n\to \infty.
\end{align*}
\begin{proof}This proof is the same analogue to the proof of Theorem
 \ref{koukaithm}. Let $f:X\to N$ be an arbitrary $1$-Lipschitz map. We
 identify the tangent space of $N$ at $z:=b\big(f_{\ast}(\mu_X)\big)$ with
 $\mathbb{R}^n$ and consider a map $f_0:=\exp_z^{-1} \circ f :X\to
 \mathbb{R}^n$. According to the hinge theorem, $f_0$ is a $1$-Lipschitz map. By using
 Lemma \ref{atatata}, we have $f_0^{-1}\big(B_{\mathbb{R}^n}\big( b\big(f_{0
 \ast}(\mu_X)\big)
 ,\rho\big)\big)=f_0^{-1}\big(B_{\mathbb{R}^n}(0,\rho)\big) =f^{-1}\big(
 B_N(z, \rho)\big)$ for any
 $\rho>0$. Hence, we obtain 
\begin{align*}
\crad(f_{\ast}(\mu_X),m-\kappa)= \crad(f_{0
 \ast}(\mu_X),m-\kappa) \leq \obsc_{\mathbb{R}^n}(X;-\kappa)
\end{align*}for any $\kappa >0$ and this completes the proof. 
\end{proof}
\end{lem}

\begin{prop}\label{shishi}Let $\{  X_n  \}_{n=1}^{\infty}$ be a sequence of mm-spaces
 of finite diameters and $p\geq 1$. We assume that a sequence of natural numbers $\{
a(n) \}_{n=1}^{\infty}$ satisfies that
\begin{align*}
\diam \Big(X_n \oblip \mathbb{R}^{a(n)},m_n-\frac{\kappa}{(\diam X_n)^p}\Big) \to 0 \text{ as }n\to \infty
\end{align*}for any $\kappa >0$. Then, for any $\kappa >0$ we have
\begin{align}\label{shinya1}
\sup \Big\{ \obsc_N\Big(X_n;-\frac{\kappa}{(\diam X_n)^p}\Big) \mid N \in \mathcal{NM}^{a(n)}      \Big\} \to
 0 \text{ as }n\to \infty.
\end{align}
\begin{proof}
Lemma \ref{eminobaka} together with Lemma \ref{morimadoka} implies that
 \begin{align}\label{ayaaya}
 & \obsc_N\Big(X_n;- \frac{\kappa}{(\diam X_n)^p}\Big)\\ \leq \ &\diam \Big(X_n
  \oblip \mathbb{R}^n,m_n-\frac{\kappa}{(\diam X_n)^p}\Big)+
  \frac{\kappa}{m_n(\diam X_n)^{p-1}} \tag*{}
  \end{align}for any $\kappa >0$. Since $p\geq 1$, (\ref{ayaaya}) leads
 to (\ref{shinya1}). This completes the proof.
\end{proof}
\end{prop}

\begin{rem}\label{redo1}\upshape
Combining Lemma \ref{kouta1} with Corollary \ref{waruiyokan} and
 Proposition \ref{shishi}, we get a simple proof of Theorem
 \ref{kameda1} in the case of $p\geq 1$. Such a way of
 proof can not be applied to obtain Theorem \ref{kameda1} in the case of
 $0<p<1$. This is because of that if $\diam X_n \to
 +\infty$ as $n\to \infty$, then the right-hand side of the inequality
 (\ref{ayaaya}) diverges to infinity for
 $0<p<1$. 
\end{rem}

\subsection{Concentration into CAT(0)-spaces and Alexandrov spaces of
  nonnegative curvature}

  The following two characterizations of CAT(0)-spaces and Alexandrov
spaces of nonnegative curvatures are due to K-T. Sturm. 
\begin{thm}[{cf.~\cite[Theorem $4.9$]{sturm}}]\label{togashi}
A complete metric space $X$ is a CAT(0)-space if and only if, for any $\nu \in
 B(X)$ with $m=
 \nu(X)$, we have
\begin{align}\label{damadama1}
 \inf_{x\in X} \int_X \dist_X(x,y)^2 \ d\nu(y)\leq \frac{1}{2m} \int\int_{X\times X} \dist_X(x,y)^2 \ d\nu(x) d\nu(y).
\end{align}
\end{thm}

\begin{thm}[{cf.~\cite[Theorem $1.4$]{sturm2}}]\label{togachan}A geodesic metric space $X$ is an
 Alexandrov space of nonnegative curvture if and only if, for $\nu\in
 B(X)$ with $m=\nu(X)$, we have
\begin{align*}
\inf_{x\in X}\int_X \dist_X(x,y)^2 \ d\nu(y)\geq \frac{1}{2m} \int\int_{X\times X} \dist_X(x,y)^2 \ d\nu(x) d\nu(y).
\end{align*}
\end{thm}

Let $X$ be a geodesic metric space. A function $\varphi : X\to \mathbb{R}$ is
called \emph{convex} if the function $\varphi \circ \gamma :[0,1]\to
\mathbb{R}$ is convex for each geodesic $\gamma :[0,1]\to X$. 
\begin{thm}[{Jensen's inequality, cf.~\cite[Theorem 6.2]{sturm}}]\label{jensen}Let $N$
 be a CAT(0)-space. Then, for any
 lower semicontinuous convex function $\varphi : N\to \mathbb{R}$ and
 any $\nu\in B^1(N)$ with $m=\nu(N)$, we have
 \begin{align*}
  \varphi (b(\nu))\leq \frac{1}{m}\int_N \varphi (x) \ d\nu (x),
  \end{align*}provided the right-hand side is well-defined.
 \end{thm}

 Since the function $\varphi:N \to \mathbb{R}$ defined by $\varphi (x):=
 \dist_N\big(x,y)^p$ for some $p\geq 1$ and $y\in N$ is the convex
 function, from Theorem \ref{jensen}, we obtain the following corollary:
 \begin{cor}For any $p\geq 1$, we have
  \begin{align}\label{damadama}
    \int_N\dist_N(b(\nu),y)^p \ d\nu(y)\leq \frac{1}{m}\int_N\int_N \dist_N(x,y)^p \ d\nu(x)d\nu(y).
   \end{align}
 \end{cor}

 Note that the above inequality (\ref{damadama}) implies worser estimate
 than the inequality (\ref{damadama1}) in the case of $p=2$. 

Let $(\Omega, \mathcal{A},\mathbb{P})$ be a probability space and $X$ a
barycentric metric space. For $X$-valued random variable $W:\Omega \to
X$ satisfying $W_{\ast}\mathbb{P} \in B^1(X)$, we define its \emph{expectation}
by the barycenter of $W_{\ast}\mathbb{P}$.

\begin{prop}\label{hanakitasaiaku}Let $N$ be a CAT(0)-space and $X$ an
 mm-space. Then for any $p\geq 1$, $\kappa >0$, and Borel measurable map
 $f:X\to N$ with $f_{\ast}(\mu_X) \in B^1(N)$, we have
 \begin{align}\label{hosiimo1}
  \crad (f_{\ast}(\mu_X),m-\kappa) \leq \frac{V_p(f)}{(m\kappa)^{1/p}}.
  \end{align}In the case of $p=2$, we also have the better estimate
\begin{align}\label{hosiimo2}
\crad (f_{\ast}(\mu_X),m-\kappa) \leq \frac{V_2(f)}{\sqrt{2m\kappa}}.
\end{align}
\begin{proof}Assume that $f_{\ast}(\mu_X)\big(N\setminus
 B_N\big(b(f_{\ast}(\mu_X)),\rho_0      \big)\big)>\kappa$ holds for $\rho_0 :=
 V_p(f)/(m\kappa)^{1/p}$. Combining (\ref{damadama}) with the Chebyshev's inequality,
 we get
\begin{align*}
\frac{V_p(f)^p}{m}=\rho_0^p \kappa <\int_N \dist_N \big(b(f_{\ast}(\mu_X)),y\big)^p \ d
 f_{\ast}(\mu_X)(y) \leq \frac{V_p(f)^p}{m},
\end{align*}which is a contradiction. Hence, we obtain $f_{\ast}(\mu_X)\big(
 B_N(b(f_{\ast}(\mu_X))  ,\rho_0       )\big)\geq m-\kappa$. This
 completes the proof of (\ref{hosiimo1}). The inequality
 (\ref{hosiimo2}) also follows from (\ref{damadama1}).
\end{proof}
\end{prop}

\begin{cor}\label{4manen}The $L^p$-concentration of maps into CAT(0)-spaces for $p\geq 1$ implies
 the concentration of maps to their expectations.
\end{cor}

From Theorem \ref{togachan}, we obtain the following corollary:
\begin{cor}Let $\{ X_n\}_{n=1}^{\infty}$ be a sequence of mm-spaces with $\sup\limits_{n\in \mathbb{N}}m_n <+\infty$ and
 $\{ Y_n\}_{n=1}^{\infty}$ a sequence of Alexandrov spaces of nonnegative
 curvatures. Assume that a sequence $\{ f_n:X_n \to Y_n
 \}_{n=1}^{\infty}$ of Borel measurable maps satisfies that
\begin{align*}
\inf_{y\in Y_n} \int_{Y_n}\dist_{Y_n}\big( y,f_n(x)   \big)^2 \ d\mu_n(x) \to 0
 \text{ as }n\to \infty.
\end{align*}Then we have $V_2(f_n) \to 0$ as $n\to \infty$.
\end{cor}

\subsection{Cases of the Gaussian concentration and the exponential concentration}

In this subsection, we review a result in \cite{ledole}. Applying
their method, we compute $\obsc_N(X;-\kappa)$ and $\obin_{N}(X)$ for
a nonpositively curved manifold
$N$ and some special mm-spaces $X$. Throughout this subsection, an mm-space $X$ always be assumed to
have a probability measure $\mu_X$.

Recall that an mm-space $(X,\dist_X,\mu_X)$ has a \emph{Gaussian concentration} whenever there are
constants $\gcon{X} >0$ and $\gconn{X}>0$ such that
\begin{align}\label{sledole1}
 \alpha_X(r)\leq \gcon{X} e^{-\gconn{X}r^2}
 \end{align}for any $r>0$. An mm-space has an \emph{exponential
 concentration} if there exist constants $\econ{X}>0$ and $\econn{X}>0$
 such that
 \begin{align}\label{sledole2}
  \alpha_X(r)\leq \econ{X} e^{- \econn{X}r}
  \end{align}for any $r>0$.

		Let $M$ be a compact connected Riemannian manifold and $\mu_M$ the normalized volume
	 measure on $M$. We shall consider $M$ as an mm-space
	 $(M,\dist_M, \mu_M)$.
\begin{thm}[{L\'{e}vy-Gromov, cf.~\cite[Section 1.2, Remark2]{milgro}
 and \cite[Theorem 2.4]{ledoux}}]\label{levygro1}Let $M$ be a compact
 connected Riemannian manifold such that $\ric_M \geq
 \widetilde{\kappa}_1>0$. Then, for any $r>0$, we have
 \begin{align*}
  \alpha_M(r)\leq e^{-\widetilde{\kappa}_1r^2/2}.
  \end{align*}
 \end{thm}
We denote by $\lambda_1(M)$ the first
	 non-zero eigenvalue of the Laplacian on $M$. 
 \begin{thm}[{Gromov-Milman, cf.~\cite[Theorem 4.1]{milgro} and \cite[Theorem 3.1]{ledoux}}]\label{gromil100}Let $M$ be a compact
 connected Riemannian manifold. Then, for any $r>0$, we have
  \begin{align*}
  \alpha_M(r)\leq e^{-\sqrt{\lambda_1(M)}r/3}.
  \end{align*}
  \end{thm}

We shall estimate $\obin_N(X)$ from the above for $N\in \mathcal{NM}^n$
and an mm-space $X$ with the Gaussian
concentration (\ref{sledole1}) or the exponential concentration (\ref{sledole2}). For the reader's
       convenience, we extract from Ledoux and Oleskiewicsz's paper
       \cite[Theorem 1]{ledole} their argument:

       Assume that an mm-space $X$ has the Gaussian concentration
       (\ref{sledole1}). From the assumption (\ref{sledole1}) and
       \cite[Proposition 1.8]{ledoux}, there exists a universal constant $C>0$, and a
      constant $C_1$ depending on $\gcon{X}$ such that
       \begin{align}\label{machigaika1}
       \mu_X( \{x\in X\mid |\varphi (x)|\geq r\}) \leq C_1 e^{-C \gconn{X}r^2}
       \end{align}for any $1$-Lipschitz function $\varphi:X\to
       \mathbb{R}$ with mean zero and $r>0$.
       \begin{rem}\upshape In the proof of \cite[Theorem 1]{ledole}, by
        citing \cite[Proposition 1.8]{ledoux}, Ledoux and Oleskiewicsz stated that the above constant $C_1$ can be chosen as the form $ C \gcon{X}$ for some
    universal constant $C >0$. But the author does not know whether
    the constant $C_1$ can be written by the above form or not only from \cite[Proposition 1.8]{ledoux}.
        \end{rem}
        We denote by $\gamma_n$ the standard
  Gaussian measure on $\mathbb{R}^n$ with density
  $(2\pi)^{-n/2}e^{-|x|^2/2}$.  For any $q\geq 0$, we put
      \begin{align*}
       M_q:= \int_{\mathbb{R}}|s|^q \ d\gamma_1 (s) = 2^{q/2}\pi^{-1/2}\Gamma\Big(\frac{q+1}{2}\Big),
       \end{align*}where $\Gamma$ is the Gamma function. By
       (\ref{machigaika1}), we get
       \begin{align}\label{riemannkk}
        \int_X |\varphi(x)|^p \ d\mu_X(x) = \ &\int_{0}^{+\infty} \mu_X( \{
        x\in X \mid |\varphi(x)|\geq r\}) \ d (r^p)\\
        \leq \ &C_1
        \int_{0}^{\infty}e^{-C \gconn{X}r^2} \ d(r^p)=
        pM_{p-1}C_1\sqrt{\frac{\pi }{2(2C \gconn{X})^{p}}} \tag*{}.
        \end{align}

        Let $f:X\to \mathbb{R}^n$ be an arbitrary $1$-Lipschitz map with
        mean zero. Observe that for every
        $y\in \mathbb{R}^n$ the map
        $y\cdot f:X\to \mathbb{R}$ is the $|y|$-Lipschitz function with
        mean zero. Hence,
        by using (\ref{riemannkk}), we have
        \begin{align}\label{riemannkk1}
         V_p(f)^p\leq \ &2^p\int_X|f(x)|^p d\mu_X(x)  \\ = \ &
         2^p\int_{X}\Big\{\frac{1}{M_p}\int_{\mathbb{R}^n}|y\cdot f(x)|^p \
         d\gamma_n(y)\Big\}d\mu_X(x) \tag*{} \\
         \leq \ &p2^pC_1\sqrt{\frac{\pi }{2(2C
         \gconn{X})^{p}}}\int_{\mathbb{R}^n} |y|^pd\gamma_n(y) \tag*{}.
         \end{align}Stirling's formula implies that 
         \begin{align*}
          \int_{\mathbb{R}^n}|y|^p \ d\gamma_n(y)=
          \frac{\sqrt{2^{p}}\Gamma(\frac{n+p}{2})}{\Gamma(\frac{n}{2})}
          \leq C (n+p)^{p/2}
          \end{align*}for some universal constant $C>0$. Therefore,
          combining this with (\ref{riemannkk1}) and Theorem \ref{koukaithm}, we obtain
          the following proposition:

 \begin{prop}\label{ledoleprop1}Let $X$ be an mm-space with the Gaussian concentration
     (\ref{sledole1}). Then, there exist a constant $C>0$ depending only on $\gcon{X}$ such that
     \begin{align*}
      \obin_N(X)  \leq C \sqrt{\frac{n+p}{\gconn{X}}}
      \end{align*}for any $n\in \mathbb{N}$, $p>0$, and $N \in
  \mathcal{NM}^n$. 
     \end{prop}

     Combining Proposition \ref{hanakitasaiaku} with Theorem \ref{levygro1} and Proposition \ref{ledoleprop1}, we get the following corollary:
     \begin{cor}\label{3manen}Let $M$ be a compact connected Riemannian manifold with
      $\ric_M\geq \widetilde{\kappa} >0$. Then, there exists a universal constant $C>0$ such that
      \begin{align*}
       \obin_N(M) \leq C\sqrt{\frac{n+p}{\widetilde{\kappa}}}
       \end{align*}for any $n\in \mathbb{N}$, $p>0$, and $N\in
      \mathcal{NM}^n$. In particular, for any $p\geq 1$ and $\kappa >0$, we have
      \begin{align}\label{orion1}
       \obsc_N(M;-\kappa) \leq \frac{C}{\kappa^{1/p}}\sqrt{ \frac{n+p}{\widetilde{\kappa}}}
       \end{align}
      \end{cor}

\begin{cor}\label{saisentan1}Let $s$ be a number with $0\leq s<1/2$ and $\{
 a(n)\}_{n=1}^{\infty}$ a sequence of natural
 numbers satisfy that $a(n)/n^{1-2s}\to 0$ as $n\to \infty$. Then for
 any $p>0$ and $\kappa >0$, we have
 \begin{align*}
  \sup \{ \obin_N(\mathbb{S}^n(n^s)) \mid N\in \mathcal{NM}^{a(n)}\} \to
  0 \text{ as }n\to \infty.
  \end{align*}In particular, for any $\kappa>0$, we have
	 \begin{align}\label{2manen}
	  \sup\{\obsc_N (\mathbb{S}^n(n^s);-\kappa) \mid N\in \mathcal{NM}^{a(n)}\}\to
	  0 \text{ as }n\to \infty.
	  \end{align}
 \begin{proof}Since $\ric_{\mathbb{S}^n(n^s)}=(n-1)n^{-2s}$, by virtue
  of Corollary \ref{3manen}, we get the corollary.
  \end{proof}
 \end{cor}

\begin{rem}\upshape Let us consider the case of $s=1/2$. Denote by $f_n$ the projection from $\mathbb{S}^n(\sqrt{n})$ into
	 $\mathbb{R}$. It is known by H. Poincar\'{e} that the sequence $\big\{   (f_{n})_{\ast} \big(\mu_{\mathbb{S}^n(\sqrt{n})}\big)    \big\}_{n=1}^{\infty}$ of probability
	 measures on $\mathbb{R}$ converges weakly to the canonical Gaussian measure on
	 $\mathbb{R}$ (cf. \cite[Section 1.1]{ledouxtalagrand}). Therefore, for any $\kappa>0$ 
	 we have
	 \begin{align*}
	  \liminf_{n \to \infty} \diam (\mathbb{S}^n (\sqrt{n})\oblip \mathbb{R},1-\kappa)>0.
	  \end{align*}
	 \end{rem}
 
 \begin{cor}\label{saisentan2}Let $\{a(n)\}_{n=1}^{\infty}$ be a sequence of natural
  numbers satisfies that $a(n)/n \to 0$ as $n\to \infty$. Then, for any
  $p>0$ and $\kappa >0$, we have
  \begin{align*}
   \sup \{  \obin_N(SO(n)) \mid N \in \mathcal{NM}^{a(n)}\} \to 0 \text{
   as }n\to \infty.
   \end{align*}
  In particular, for any $\kappa >0$, we have
  \begin{align}\label{7manen}
   \sup \{  \obsc_N(SO(n);-\kappa) \mid N\in \mathcal{NM}^{a(n)}\} \to
   0\text{ as }n\to \infty.
   \end{align}
  \begin{proof}Since $\ric_{SO(n)}\geq (n-1)/4$, we obtain the corollary.
   \end{proof}
  \end{cor}

  In \cite[Theorem 1]{ledole}, by using (\ref{riemannkk1}), Ledoux and Oleszkievwicz
  proved the following theorem:

  \begin{thm}[{Ledoux-Oleszkievwicz, cf.~\cite[Theorem 1]{ledole}}]\label{thledole}Let $X$ be an mm-space with
   the Gaussian concentration (\ref{sledole1}). Then, there exists a
   universal constant $C>0$ and a constant $C_1>0$ depending only on $\gcon{X}>0$ such that
   \begin{align*}
    f_{\ast}(\mu_X) \big( N\setminus B_N \big( b(f_{\ast}(\mu_X)),r\big)\big) \leq C_1 \gamma_n (\{ x\in \mathbb{R}^n
    \mid |x|\geq C \sqrt{\gconn{X}}r\})
    \end{align*}for any $n \in \mathbb{N}$, $r>0$, $N\in \mathcal{NM}^n$,
   and
   $1$-Lipschitz map $f:X\to N$.
   \end{thm}
   In \cite{ledole}, Ledoux and Oleszkievwicz stated Theorem
    \ref{thledole} in the situation of $N=\mathbb{R}^n$. The general
    case follows from the same way of the proof of Lemma
    \ref{morimadoka}.

    Let $\{ X_n\}_{n=1}^{\infty}$ be a sequence of mm-spaces having the
      Gaussian concentration (\ref{sledole1}). Assume that
      $\sup\limits_{n\in \mathbb{N}} \gcon{X_n} <+\infty$ and a
      sequence $\{
      a(n)\}_{n=1}^{\infty}$ of natural numbers satisfies that $a(n)/
      \gconn{X_n} \to 0$ as $n\to \infty$. In this situation, by using
      Corollary \ref{4manen} and Proposition \ref{ledoleprop1}, we get
      \begin{align}\label{realkodaemotion}
       \sup \{  \obsc_N(X_n;-\kappa) \mid N\in \mathcal{NM}^{a(n)}\} \to
       0 \text{ as }n\to \infty.
       \end{align}We note that (\ref{realkodaemotion}) also follows from Theorem
       \ref{thledole} and
       the following lemma:
    \begin{lem}[{cf.~\cite[Corollary 2.3]{ale}}]For any $r\geq \sqrt{n}$, we have 
     \begin{align*}
      \gamma_n ( \{  x\in \mathbb{R}^n \mid |x|\geq r\})\leq \exp \Big(
      -\frac{n}{4} \Big(1-\frac{n}{r^2}\Big)^2\Big).
      \end{align*}
     \end{lem}

    Assume that an mm-space has the exponential concentration
    (\ref{sledole2}). From the assumption (\ref{sledole2}) and
    \cite[Proposition 1.8]{ledoux}, there exist a universal constant
    $C>0$, and a constant $C_1$ depending on $\econ{X}$ such that
    \begin{align*}
     \mu_X (\{  x\in X \mid |\varphi (x)|\geq r\})\leq C_1 e^{-C \econn{X}r }
     \end{align*}for any $1$-Lipschitz function $\varphi :X\to
     \mathbb{R}$ with mean zero and $r>0$. Hence, by the same method in the
     proof of Proposition \ref{ledoleprop1}, we obtain the following proposition:

     \begin{prop}\label{exnoprop}Let $X$ be an mm-space with the exponential concentration
     (\ref{sledole2}). Then, there exist a constant $C>0$ depending only on $\econ{X}$ such that
     \begin{align*}
      \obin_N(X)  \leq C \frac{\sqrt{p(n+p)}}{\econn{X}}
      \end{align*}for any $n\in \mathbb{N}$, $p\geq 1$, and $N\in \mathcal{NM}^n$.
      \end{prop}

      Proposition \ref{hanakitasaiaku} together with Theorem \ref{gromil100}
      and Proposition \ref{exnoprop} implies the
      following corollary:
      \begin{cor}\label{lamdadada}Let $M$ be a compact connected Riemannian
       manifold. Then, there exists a universal constant $C>0$ such that
       \begin{align}\label{10manen}
        \obin_N(M)\leq C \sqrt{\frac{p(n+p)}{\lambda_1(M)}}
        \end{align}for any $n\in \mathbb{N}$, $p\geq 1$, and $N\in
       \mathcal{NM}^n$. In particular, for any $\kappa >0$, we have
       \begin{align}\label{orion2}
        \obsc_N(M;-\kappa)\leq \frac{C}{\kappa^{1/p}}\sqrt{\frac{p(n+p)}{\lambda_1(M)}}.
        \end{align}
      \end{cor}

      \begin{cor}\label{lamdadada1}If a sequence $\{
      M_n\}_{n=1}^{\infty}$ of compact connected Riemannian manifolds
      with and a sequence $\{
      a(n)\}_{n=1}^{\infty}$ of natural numbers satisfy that
      $a(n)/\lambda_1(M_n) \to 0$ as $n\to \infty$, then for any
      $p>0$ we have
      \begin{align*}
       \sup \{ \obin_N(M_n) \mid N\in \mathcal{NM}^{a(n)}\} \to 0 \text{
       as }n\to \infty.
       \end{align*}In particular, for any $\kappa >0$, we have
       \begin{align*}
        \sup \{\obsc_N (M_n;-\kappa) \mid N\in \mathcal{NM}^{a(n)}\} \to
        0 \text{ as }n\to \infty.
        \end{align*}
       \end{cor}
      Applying the Lichnerowicz theorem of the first eigenvalue of
      the Laplacian to
      Corollary \ref{lamdadada1}, we get the same conclusions of
      Corollaries \ref{saisentan1} and \ref{saisentan2}.

\subsection{Gromov's results}In this subsection, we review some Gromov's
      results.

      \begin{thm}[{Gromov's isoperimetry of waists, cf.~\cite{gromov2}}]\label{isowa}Let $m$ and $n$ are natural
       numbers with $m\geq n$ and $f:\mathbb{S}^m \to
       \mathbb{R}^n$ be a continuous map. Then, there
 exists a point $m_f \in \mathbb{R}^n$ such that 
       \begin{align*}
        \mu_n \big( \big(f^{-1}(m_f)\big)_{\varepsilon}\big)\geq \mu_n \big((\mathbb{S}^{m-n})_{\varepsilon}\big)
        \end{align*}for any $\varepsilon >0$, where
       $\mathbb{S}^{m-n}\subseteq \mathbb{S}^m$ denotes an equatorial (m-n)-sphere.
       \end{thm}

       By using Theorem \ref{isowa}, one might be able to examine the
       asymptotic behavior of the observable diameters $\diam
       (\mathbb{S}^n \oblip \mathbb{R}^{a(n)},1-\kappa)$.
\begin{thm}[{Gromov, cf.~\cite[Section $8$]{gromovcat}}]\label{sugeyo}Let $M$ be an $m$-dimensional compact symmetric space of rank
 $1$ and $N$ be an $n$-dimensional Riemannian manifold. Then, for any $1$-Lipschitz map
 $f:M\to N$, we have
\begin{align*}
V_2(f)^2\leq \frac{n}{m}  \int \int_{M\times M}
 \dist_M (x,x')^2 \ d\mu_M(x) d\mu_M(x') .
\end{align*}
\end{thm}
We indicate by $\mathcal{M}^n$ the set of all $n$-dimensional Riemannian
manifolds.
\begin{cor}\label{kouta4}Assume that $0\leq s< 1/2$ and a sequence $\{  a(n)
 \}_{n=1}^{\infty}$ of natural numbers satisfies $a(n)/n^{1-2s}\to 0$ as
 $n\to \infty$. Then we have
\begin{align}\label{2107}
\sup \{  \obinin_{M}\big(\mathbb{S}^n(n^s)\big) \mid M \in
 \mathcal{M}^{a(n)}           \} \to 0 \ \text{as } n\to \infty.
\end{align}In particular, for any $\kappa>0$ we have
\begin{align*}
\sup \{  \diam(\mathbb{S}^n(n^s) \oblip M,1 -\kappa) \mid M \in
 \mathcal{M}^{a(n)}           \} \to 0 \ \text{as } n\to \infty.
\end{align*}
\end{cor}

\begin{rem}\upshape Since compact symmetric spaces $SO(n)$ are not rank $1$ for $n\geq 4$, we
 cannot apply Theorem \ref{sugeyo} for $SO(n)$.
\end{rem}

Combining (\ref{2107}) with Corollary \ref{4manen}, we also get (\ref{2manen}).

	\begin{thm}[{Gromov, cf.~\cite[Section
	 $3\frac{1}{2}.41$]{gromov}}]\label{g}Let $M$ be a compact
	 connected Riemannian manifold and $N\in \mathcal{NM}^n$. Then, for any $\kappa>0$ we have
	 \begin{align}\label{215}
	  \obsc_N(M;-\kappa)\leq \frac{1}{\sqrt{\kappa}}\sqrt{\frac{n}{\lambda_1(M)}}.
	  \end{align}
\begin{proof}In \cite[Section $3\frac{1}{2}.41$]{gromov}, Gromov proved
 the case of $N=\mathbb{R}^n$. This and Lemma \ref{morimadoka} gives the proof
 of the theorem.
\end{proof}
	 \end{thm}From Theorem \ref{g}, we obtain (\ref{2manen}) and (\ref{7manen}). Compare the above inequality (\ref{215}) with (\ref{orion1}) and (\ref{orion2}).
\subsection{Applications of Theorem \ref{th}}
The proof of the following lemma is easy.
	\begin{lem}[{cf.~\cite[Section $3\frac{1}{2}.32$]{gromov}}]\label{yaku}For an mm-space $X $, $n\in \mathbb{N}$, and $\kappa >0$, we have
	 \begin{align*}
	  \diam (X\oblip \mathbb{R}^n,m-n\kappa)\leq \sqrt{n}\diam (X\oblip \mathbb{R},m-\kappa).
	  \end{align*}
	 \end{lem}

\begin{lem}\label{damasenai}Let $M$ be a compact connected Riemannian
 manifold such that $Ric_M \geq \widetilde{\kappa}_1 >0$. Then, for any
 $\kappa>0$ we have
\begin{align*}
\diam (M \oblip \mathbb{R}^n,1-\kappa)\leq 2\sqrt{\frac{2n\log\big(\frac{2n}{\kappa}\big)}{\widetilde{\kappa}_1}}.
\end{align*}
\begin{proof}Theorem \ref{levygro1} together with Lemmas \ref{ccc}, \ref{a1}, and
 \ref{yaku} leads to the proof of the lemma.
\end{proof}
\end{lem}

\begin{rem}\upshape For fixed $\kappa >0$, the inequality in Lemma
 $\ref{damasenai}$ is weaker than that the inequalities (\ref{orion1}),
 (\ref{orion2}), and (\ref{215}) for a
 Riemannian manifold $M$ with a positive lower Ricci curvature
 bound for the same screen $\mathbb{R}^n$. However, in the case where
 $\kappa \to 0$, the inequality in Lemma
 $\ref{damasenai}$ is shaper than (\ref{orion1}), (\ref{orion2}), and
 (\ref{215}). 
\end{rem}

\begin{cor}Let $0\leq s <1/2$. Assume that sequences $\{
	 a(n)\}_{n=1}^{\infty} $, $\{
	 \kappa_n\}_{n=1}^{\infty}$, $\{s_n\}_{n=1}^{\infty}$ of real
 numbers satisfy that $a(n)\in \mathbb{N}$, $s_n\geq 0$, $\kappa_n<0$, $\sup\limits_{n\in \mathbb{N}}s_n
  <1/2$, $\sup\limits_{n\in \mathbb{N}}\kappa_n <0$,
\begin{align}\label{sugakutou5}
 \frac{a(n)\log n}{n^{1-2s}} \to 0 \text{ as }n\to \infty, \text{ and }
  \frac{a(n)\log (-\kappa_n)^{s_n}}{n^{1-2s}}\to 0 \text{ as }n\to \infty.
 \end{align}Then for any 
	 sequences $\{N_n\}_{n=1}^{\infty}$ with $N_n \in
	      \mathcal{NM}^{a(n)}(s_n;\kappa_n)$ and $\{ f_n:X_n \to
	      N_n\}_{n=1}^{\infty} $ of $1$-Lipschitz maps, we have
	      \begin{align*}
	       \diam I_{N_n}\big(  (f_n)_{\ast}(\mu_{\mathbb{S}^n(n^s)})
	       ;1-2s_n,1     \big) \to 0 \text{ as }n\to \infty.
	       \end{align*}
\begin{proof}
Since $\ric_{S^n(n^s)}=(n-1)n^{-2s}$, applying Lemma \ref{damasenai},
for any $\kappa>0$  we obtain
\begin{align*}
\diam \Big(\mathbb{S}^n(n^s) \oblip \mathbb{R}^{a(n)},1-\frac{\kappa}{(-\kappa_n)^{s_n}\diam
 \mathbb{S}^n(n^s)} \Big) \leq 2\sqrt{\frac{2a(n)\log \Big(\frac{2a(n)(-\kappa_n)^{s_n} n^{s}\pi}{\kappa}\Big)}{(n-1)n^{-2s}}                     }.
\end{align*}Therefore, from Theorem \ref{th}, this completes the proof. 
\end{proof}
	     \end{cor}

         \begin{cor}Let $0\leq s <1/3$. Assume that sequences $\{
	 a(n)\}_{n=1}^{\infty} $, $\{
	 \kappa_n\}_{n=1}^{\infty}$, $\{s_n\}_{n=1}^{\infty}$ of real
 numbers satisfy that $a(n)\in \mathbb{N}$, $s_n\geq 0$, $\kappa_n<0$, $\sup\limits_{n\in \mathbb{N}}s_n
  <1/2$, $\sup\limits_{n\in \mathbb{N}}\kappa_n <0$,
          \begin{align}\label{sugakutou6}
           \frac{a(n)}{n^{1-2s}}\to 0 \text{ as }n\to \infty
           \text{,
           and } \frac{(-\kappa_n)^{s_n}}{n^{1-3s}}\to 0 \text{ as }n\to
           \infty. 
           \end{align}Then for any 
	 sequences $\{N_n\}_{n=1}^{\infty}$ with $N_n \in
	      \mathcal{NM}^{a(n)}(s_n;\kappa_n)$ and $\{ f_n:X_n \to
	      N_n\}_{n=1}^{\infty} $ of $1$-Lipschitz maps, we have
	      \begin{align*}
	       \diam I_{N_n}\big(  (f_n)_{\ast}(\mu_{\mathbb{S}^n(n^s)})
	       ;1-2s_n,1     \big) \to 0 \text{ as }n\to \infty.
	       \end{align*}
          \begin{proof}We may assume that $\kappa_n \leq -1$ for any
	       $n\in \mathbb{N}$. Since $p_n:=(-\kappa_n)^{s_n} \diam
	       \mathbb{S}^n(n^s)=n^s (-\kappa_n)^{s_n}\pi \geq 1$, combining
	       Lemma \ref{hayakusitekure} with (\ref{orion1}), we have 
           \begin{align*}
            \diam \Big( \mathbb{S}^n(n^s) \oblip
            \mathbb{R}^{a(n)},1-\frac{\kappa}{(-\kappa_n)^{s_n}\diam
            \mathbb{S}^n (n^s)}\Big) \leq
            \frac{Cp_n^{1/p_n}}{\kappa^{1/p_n}}\sqrt{ \frac{a(n)+ p_n}{(n-1)n^{-2s}}}
            \end{align*}for any $\kappa>0$. Hence, from Theorem \ref{th}, this completes
           the proof.
           \end{proof}
          \end{cor}
          \begin{rem}\upshape Compare the assumption (\ref{sugakutou5}) with (\ref{sugakutou6}). One can take
           $a(n)$ in (\ref{sugakutou6}) to be greater than in
           (\ref{sugakutou5}), whereas one can not take $(-\kappa_n)$ in
           (\ref{sugakutou6}) to be greater than $(-\kappa_n)$ in
           (\ref{sugakutou5}).
           \end{rem}
\begin{cor}Assume that sequences $\{
	 a(n)\}_{n=1}^{\infty} $, $\{
	 \kappa_n\}_{n=1}^{\infty}$, $\{s_n\}_{n=1}^{\infty}$ of real
 numbers satisfy that $a(n)\in \mathbb{N}$, $ s_n \geq 0$,
 $\sup\limits_{n\in \mathbb{N}}s_n <1/2$, $\sup\limits_{n\in
 \mathbb{N}}\kappa_n <0$,
 \begin{align}\label{ochahara11}
  \frac{a(n)\log n}{n} \to 0 \text{ as }n\to \infty, \text{
  and }  \frac{a(n)\log(-\kappa_n)^{s_n}}{n} \to 0 \text{ as }n\to \infty.
  \end{align}Then for any 
	 sequences $\{N_n\}_{n=1}^{\infty}$ with $N_n \in
	      \mathcal{NM}^{a(n)}(s_n;\kappa_n)$ and $\{ f_n:X_n \to
	      N_n\}_{n=1}^{\infty} $ of $1$-Lipschitz maps, we have
	      \begin{align*}
	       \diam I_{N_n} \big(  (f_n)_{\ast}(\mu_{SO(n)})
	       ;1-2s_n,1     \big)\to 0 \text{ as }n\to \infty.
	       \end{align*}
\begin{proof}
Since $\ric_{SO(n)}\geq (n-1)/4$, by the Myers's diameter theorem, we
get $\diam SO(n) \leq \pi \sqrt{2n}$. From this and Lemma
\ref{damasenai}, for any $\kappa >0$, we have 
\begin{align*}
\diam \Big(  SO(n)\oblip \mathbb{R}^{a(n)}
	   ,1-\frac{\kappa}{(-\kappa_n)^{s_n}\diam SO(n)}
	   \Big)\leq 2\sqrt{\frac{8a(n)\log \Big(  \frac{2 a(n)(-\kappa_n)^{s_n}\sqrt{2n}\pi}{\kappa}             \Big)}{n-1}}. 
\end{align*}Hence, applying Theorem \ref{th}, this completes the proof.
\end{proof}
	     \end{cor}

         \begin{cor}Assume that sequences $\{
	 a(n)\}_{n=1}^{\infty} $, $\{
	 \kappa_n\}_{n=1}^{\infty}$, $\{s_n\}_{n=1}^{\infty}$ of real
 numbers satisfy that $a(n)\in \mathbb{N}$, $ s_n \geq 0$,
 $\sup\limits_{n\in \mathbb{N}}s_n <1/2$, $\sup\limits_{n\in
 \mathbb{N}}\kappa_n <0$,
 \begin{align}\label{ochahara12}\frac{a(n)}{n} \to 0 \text{ as }n\to \infty, \text{
  and }  \frac{(-\kappa_n)^{s_n}}{\sqrt{n}} \to 0 \text{ as }n\to \infty.
  \end{align}Then for any 
	 sequences $\{N_n\}_{n=1}^{\infty}$ with $N_n \in
	      \mathcal{NM}^{a(n)}(s_n;\kappa_n)$ and $\{ f_n:X_n \to
	      N_n\}_{n=1}^{\infty} $ of $1$-Lipschitz maps, we have
	      \begin{align*}
	       \diam I_{N_n} \big(  (f_n)_{\ast}(\mu_{SO(n)})
	       ;1-2s_n,1     \big)\to 0 \text{ as }n\to \infty.
	       \end{align*}
          \begin{proof}Since $p_n:= \pi (-\kappa_n)^{s_n}\sqrt{2n} \geq 1$ for any sufficiently
	       large $n\in \mathbb{N}$, combining Lemma \ref{hayakusitekure}
	       with (\ref{orion1}), we get
           \begin{align*}
            \diam \Big(SO(n)\oblip \mathbb{R}^{a(n)},
            1-\frac{\kappa}{(-\kappa_n)^{s_n}\diam SO(n)}\Big) \leq
            \frac{Cp_n^{1/p_n}}{\kappa^{1/p_n}} \sqrt{\frac{a(n)+p_n}{n-1}}
            \end{align*}for any $\kappa >0$. Therefore, from Theorem \ref{th}, this
           completes the proof.
           \end{proof}
          \end{cor}

          Compare the assumption (\ref{ochahara11}) with (\ref{ochahara12}).

	\section{Tree screens}
	We define a \emph{tree} $T$ as a (possibly infinite) connected combinatorial graph having no loops. We identify the individual
	edges of a tree as bounded closed intervals of the real lines, and then define the
	distance between two points of the tree to be the infimum of the lengths of paths joining them. 
	\begin{dfn}\upshape  Let $(X,\dist_X,\mu_X)$ be an mm-space and $f:X\to T$ be a Borel measurable map. A \emph{pre-L\'{e}vy
	 mean} of $f$ is a
	 point $p\in T$ such that there exist two trees $T',T'' \subseteq T$ such that
	 \begin{align*}
	  T=T'\cup T'', \ T'\cap T'' =\{ p  \}, \ f_{\ast}(\mu_X)(T')\geq \frac{m}{3},\ \text{ and }f_{\ast}(\mu_X)(T'')\geq \frac{m}{3}.
	  \end{align*}
	 \end{dfn}
	\begin{prop}There exists a pre-L\'{e}vy mean.
	 \begin{proof}Take an edge $e$ of $T$ and fix an inner point $q\in e$. There exist two trees $T',T'' \subseteq T$ such that
	  $T'\cap T'' =\{ q  \}$ and $T=T' \cup T''$. If $f_{\ast}(\mu_X)(T'),f_{\ast}(\mu_X)(T'')\geq m/3$, we have finished the proof.
	  Hence we consider the case of $f_{\ast}(\mu_X)(T'')<m/3$.

	  Let $V'$ be the vertex set of $T'$. For any $v \in V'$, we indicate by $\mathcal{C}_v$ the set of all connected components
	  of $T \setminus \{ v \}$ and put $\mathcal{C}_v':=
	  \big\{ \widetilde{T} \cup \{  v   \} \mid \widetilde{T} \in \mathcal{C}_v  \big\}$. Suppose that a point $v\in V'$
	  satisfies $f_{\ast}(\mu_X)(\widetilde{T})< m/3$ for any $\widetilde{T} \in
	  \mathcal{C}_{v}'$, then it is easy to check that $v$ is a pre-L\'{e}vy mean of $f$.
	  So, we assume that for any $v\in V'$ there exists $T_v \in \mathcal{C}_v'$ such that $f_{\ast}(\mu_X)(T_v)\geq
	  m/3$. If for some $v \in V'$
	  there exists $ T_v' \in \mathcal{C}_v' \setminus \{  T_v   \}$ such that  $f_{\ast}(\mu_X)(T_v')\geq m/3$, then this 
	  $v$ is a  pre-L\'{e}vy mean of $f$. Therefore, we also assume that 
	  $f_{\ast}(\mu_X)(T_v')<m/3$ for any $v\in V'$ and $T_v' \in \mathcal{C}_v' \setminus \{  T_v   \}$. 

	  We denote by $\Gamma$ the set of all unit speed geodesics $\gamma:[0,L(\gamma)]\to T'$ from $q$ such that $\gamma \big(L (\gamma) \big) \in V'$ and $\gamma \big(  [t,L(\gamma)]\big) \subseteq T_{\gamma(t)}$ for each
	  $\gamma(t)\in\big( V'\setminus \big\{  \gamma \big(L (\gamma)\big)   \big\} \big)\cap \gamma\big([0,L(\gamma)]\big)$. 
	  It is easy to verify that $\gamma \subseteq \gamma'$ for any $\gamma ,\gamma' \in \Gamma$ with $L(\gamma)\leq
	  L(\gamma')$. Put $\alpha : = \sup \{  L(\gamma)  \mid \gamma \in \Gamma\}$. Let us show that there exists $\widetilde{\gamma} \in
	  \Gamma$ with $L(\widetilde{\gamma})= \alpha$. If $L(\gamma)<\alpha$ for any $\gamma \in \Gamma$, there exists a
	  sequence $\{  \gamma_n   \}_{n=1}^{\infty} \subseteq \Gamma$ such that $L(\gamma_1)<L (\gamma_2)< \cdots \to \alpha$ as
	  $n\to \infty$. Then we have
	  \begin{align*}
	   f_{\ast}(\mu_X)(T)=\limsup_{n \to \infty} f_{\ast}(\mu_X)\Big\{  \Big( \bigcup \mathcal{C}_{\gamma_n (L(\gamma_n))}'\Big)\setminus
	   T_{\gamma_n(L(\gamma_n))}      \Big\}\leq \frac{2m}{3},
	   \end{align*}which is a contradiction. Suppose that there exists a sequence $\{  t_n
	  \}_{n=1}^{\infty} \subseteq [0,L(\widetilde{\gamma})]$ such that $t_1 <t_2 < \cdots \to L( \widetilde{\gamma})$ as $n\to \infty$ and $\widetilde{\gamma} (t_n)
	  \in V'$ for each $n\in \mathbb{N}$. Since $  T'' \subseteq T_{\widetilde{\gamma}(L(\widetilde{\gamma}))} $ and 
	  \begin{align*}
	   f_{\ast}(\mu_X)\Big\{ \Big(\bigcup \mathcal{C}_{\widetilde{\gamma}(L(\widetilde{\gamma}))} \setminus  T_{\widetilde{\gamma}(L(\widetilde{\gamma}))} \Big) \cup \big\{ \widetilde{\gamma} \big(
	   L (\widetilde{\gamma})\big)   \big\}\Big\} = \lim_{n \to \infty} f_{\ast}(\mu_X) (T_{\widetilde{\gamma}(t_n)}) \geq \frac{m}{3},
	   \end{align*}$\widetilde{\gamma} \big(L(\widetilde{\gamma})\big)$ is a pre-L\'{e}vy mean of $f$. We will consider the other
	  case, that is, there exist $t_0 \in [0,L(\widetilde{\gamma})]$ and edge $e_0$ of $T'$ such that $\widetilde{\gamma}(t_0)\in
	  V'$ and $e_0$ connects $\widetilde{\gamma} (t_0)$ and $\widetilde{\gamma}\big(L(\widetilde{\gamma})\big)$.
	  If $f_{\ast}(\mu_X)\big( (\bigcup \mathcal{C}_{\widetilde{\gamma} (t_0)} \setminus T_{\widetilde{\gamma} (t_0)} )\cup \{ \widetilde{\gamma}(t_0) \} \big) \geq
	  m/3$, then $\widetilde{\gamma} (t_0)$ is a pre-L\'{e}vy mean of $f$. If $f_{\ast}(\mu_X)\big( (\bigcup \mathcal{C}_{\widetilde{\gamma} (t_0)} \setminus T_{\widetilde{\gamma} (t_0)} )\cup \{
	  \widetilde{\gamma}(t_0) \} \big) < m/3$, there exists a pre-L\'{e}vy
	  mean of $f$ on $e_0$ since $e_0 \subseteq T_{\widetilde{\gamma}(L(\widetilde{\gamma}))}$. This completes the proof.

	 \begin{figure}[tbp]
	 \begin{center}
	  \includegraphics[width=7cm,clip]{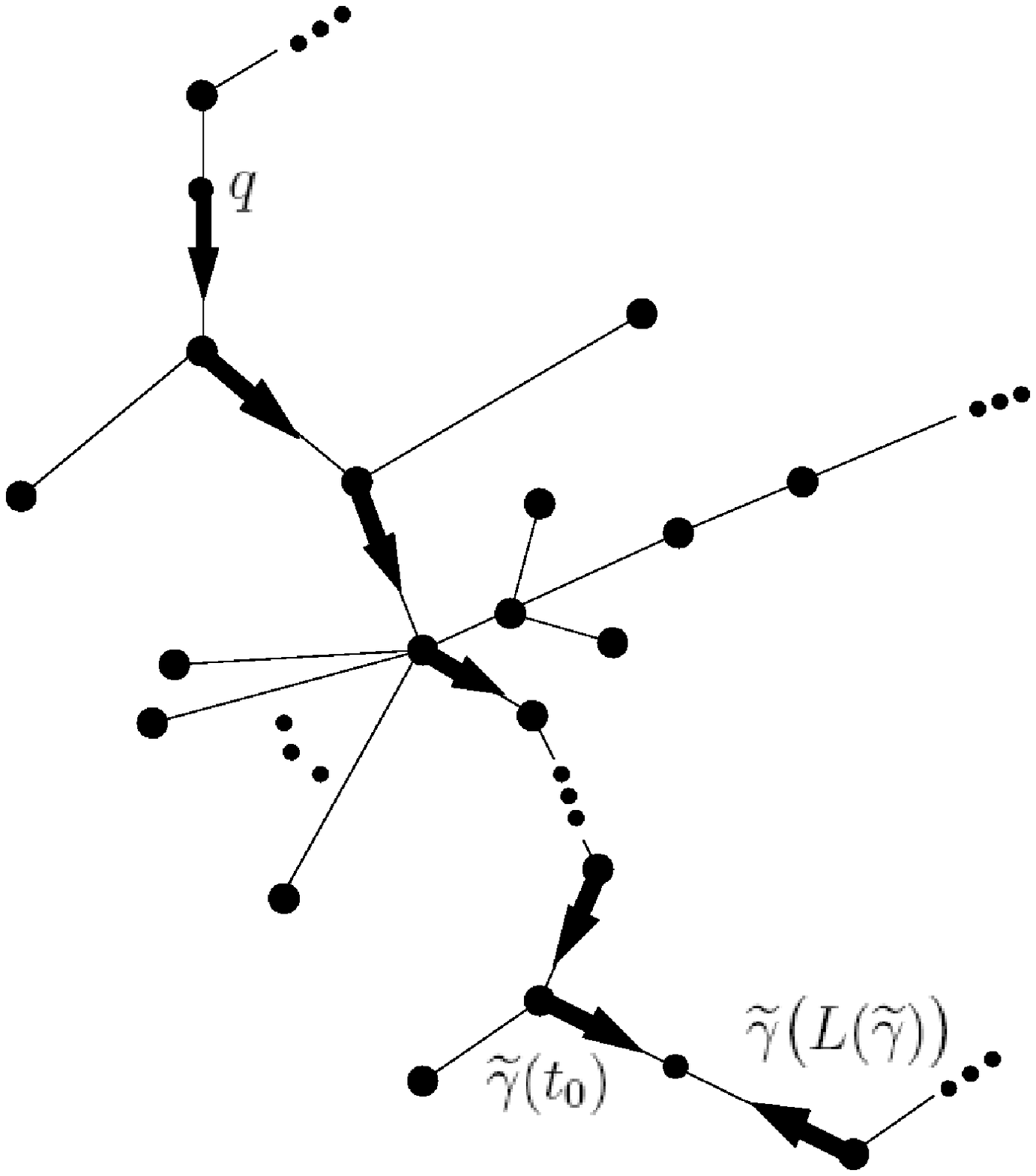}
	  \caption{Tree}
	  \label{notBG}
	 \end{center}
	  \end{figure}

	  \end{proof}
	 \end{prop}

\begin{lem}\label{ainouta}For any $\kappa>0$, we have
\begin{align*}
\diam (X \oblip T,m-\kappa) \leq 2\sep\Big( X;\frac{m}{3},\frac{\kappa}{2} \Big)\end{align*}

	 \begin{proof}Let $ f :X \to T $ be an arbitrary $1$-Lipschitz map. Take a pre-L\'{e}vy mean $p\in T$ of $f$ and let $T',T'' \subseteq T$
	  be its associated trees. Let $\varepsilon >0$ satisfies
	  $\varepsilon > \sep(\mu_X;m/3,\kappa/2)$. Suppose that $f_{\ast}(\mu_X)\big(T'\setminus
	   B_{T}(p,\varepsilon))\geq \frac{\kappa}{2}$. Then, applying Lemma \ref{neko}, we have
	  \begin{align*}
	    \varepsilon  \leq \dist_{T} (T'',T'\setminus
	   B_{T}(p,\varepsilon) )\leq \sep
	   \Big(f_{\ast}(\mu_X);\frac{m}{3},\frac{\kappa}{2} 
	   \Big)\leq \sep \Big(\mu_X; \frac{m}{3},\frac{\kappa}{2}  \Big),
	   \end{align*}which is a contradiction. 
	  In the same way, we have $f_{\ast}(\mu_X)\big(  T''\setminus B_{T}(p,\varepsilon)    \big)<\kappa/2$. As a
	  consequence, we obtain $f_{\ast}(\mu_X)\big(  T\setminus
	  B_{T}(p,\varepsilon)    \big)<\kappa$. This completes the proof.
	  \end{proof}
\end{lem}

\begin{proof}[Proof of Proposition \ref{treenotoki}]The claim obviously
 follows from Lemma \ref{ainouta}.
\end{proof}

	\begin{ack}\upshape
	 The author would like to thank Professor Takashi Shioya for his
	 valuable suggestions related to Theorem \ref{th} and many discussions. He thanks to Professor Vitali Milman
	 for useful comments. He also thanks to the referee for carefully
	 reading the manuscript and fruitful suggestions. Without them, this work would have never been completed.
	 \end{ack}

	\end{document}